\documentclass[11pt, letterpaper]{article}
\pdfoutput=1

\usepackage[authoryear]{natbib}
\usepackage{amsmath, amssymb, amsfonts, amsbsy, amsthm, epsfig, subfig, graphicx,rotating, tabularx, booktabs, setspace}

\usepackage{fancyhdr}

\RequirePackage[colorlinks,citecolor=blue,urlcolor=blue]{hyperref}

\makeatletter
\def\singlespace{\def\baselinestretch{1}\@normalsize}

\makeatletter
\def\singlespace{\def\baselinestretch{1}\@normalsize}

\numberwithin{equation}{section}

\renewcommand{\hat}{\widehat}

\newcommand{\bfm}[1]{\ensuremath{\mathbf{#1}}}

     \def\bD{\bfm D}          
\def\be{\bfm e}               
              
     \def\bG{\bfm G}          
     \def\bH{\bfm H}          
     \def\bI{\bfm I}

     \def\bM{\bfm M}

     \def\bP{\bfm P}

\def\bu{\bfm u}     \def\bU{\bfm U}          
\def\bv{\bfm v}     \def\bV{\bfm V}          
               
\def\bx{\bfm x}     \def\bX{\bfm X}          
     \def\bY{\bfm Y}          
\def\bz{\bfm z}     \def\bZ{\bfm Z}

\newcommand{\bfsym}[1]{\ensuremath{\boldsymbol{#1}}}
\def \balpha   {\bfsym{\alpha}}       \def \bbeta    {\bfsym{\beta}}
       
     \def \bzeta    {\bfsym{\zeta}}
         \def \btheta   {\bfsym{\theta}}

          \def \bxi      {\bfsym{\xi}}

\def \beps     {\bfsym \varepsilon}

\def \bSigma   {\bfsym{\Sigma}}

\renewcommand{\hat}{\widehat}

\def \heps     {\hat{\heps}}

       \def \hbbeta    {\hat{\bbeta}}

\DeclareMathOperator*{\argmin}{argmin}
\DeclareMathOperator{\corr}{corr}

\DeclareMathOperator{\sgn}{sgn}
\DeclareMathOperator{\supp}{supp}

\def \var   {{\rm Var}}
\def \corr  {\mbox{corr}}

\def \sgn   {\mbox{sgn}}

 \def\P{\mathbb{P}}

 at 8truept

\def\today{\ifcase\month\or
  January\or February\or March\or April\or May\or June\or
  July\or August\or September\or October\or November\or December\fi
  \space\number\day, \number\year}

\def \wt      {\widetilde}

\def \no      {\noindent}
\def \newpage {\vfill\eject}

\newdimen\biblioindent\biblioindent=30pt
\newcommand{\beq}  {\begin{equation}}
\newcommand{\eeq}  {\end{equation}}
\newcommand{\beqn} {\begin{eqnarray}}
\newcommand{\eeqn} {\end{eqnarray}}
\newcommand{\beqnn}{\begin{eqnarray*}}
\newcommand{\eeqnn}{\end{eqnarray*}}

\def \wt {\widetilde}

\def\bz{\bfm z}

\newcommand{\bbX}{\mathbb{X}}
\newcommand{\mo}{\mathbf{0}}
\newcommand{\sn}{\sum_{i=1}^n}
\newcommand{\sump}{\sum_{j=1}^p}
\newcommand{\e}{\mathbb{E}}
\newcommand{\argmax}{\mathop{\rm arg\max}}
\def\bbr{ \mathbb{R}}
\newcommand{\nn}{\nonumber}
\newtheorem{theorem}{Theorem}[section]
\newtheorem{corollary}{Corollary}[section]
\newtheorem{remark}{Remark}[section]
\newtheorem{assumption}{Condition}[section]
\newtheorem{lemma}{Lemma}[section]

\allowdisplaybreaks
\renewcommand{\baselinestretch}{1.66}
\newtheorem{proposition}{Proposition}[section]
\newcounter{CondCounter}

\def\T{\mathrm{\scriptstyle T}} 

\begin{document}

\title{Guarding against Spurious Discoveries in High Dimensions}

\author{Jianqing Fan\footnote{Department of Operations Research and Financial Engineering, Princeton University, Princeton, NJ 08544. The research of Jianqing Fan was supported in part by NSF Grants DMS-1206464,  DMS-1406266, and NIH Grant R01-GM072611-10.} ~~and ~Wen-Xin Zhou\footnote{Department of Operations Research and Financial Engineering, Princeton University, Princeton, NJ 08544.}}

\date{}

\maketitle
\begin{abstract}
Many data mining and statistical machine learning algorithms have been developed to select a subset of covariates to associate with a response variable. Spurious discoveries can easily arise in high-dimensional data analysis due to enormous possibilities of such selections. How can we know statistically our discoveries better than those by chance? In this paper, we define a measure of goodness of spurious fit, which shows how good a response variable can be fitted by an optimally selected subset of covariates under the null model, and propose a simple and effective LAMM algorithm to compute it. It coincides with the maximum spurious correlation for linear models and can be regarded as a generalized maximum spurious correlation. We derive the asymptotic distribution of such goodness of spurious fit for generalized linear models and $L_1$ regression. Such an asymptotic distribution depends on the sample size, ambient dimension, the number of variables used in the fit, and the covariance information. It can be consistently estimated by multiplier bootstrapping and used as a benchmark to guard against spurious discoveries. It can also be applied to model selection, which considers only candidate models with goodness of fits better than those by spurious fits. The theory and method are convincingly illustrated by simulated examples and an application to the binary outcomes from German Neuroblastoma Trials.
\end{abstract}


\newpage

\section{Introduction}
\label{sec:introduction}

Technological developments in science and engineering lead to collections of massive amounts of high-dimensional data. Scientific advances have become more and more data-driven, and researchers have been making efforts to understand the contemporary large-scale and complex data. Among these efforts, variable selection plays a pivotal role in high-dimensional statistical modeling, where the goal is to extract a small set of explanatory variables that are associated with given responses such as biological, clinical, and societal outcomes. Toward this end, in the past two decades, statisticians have developed many data learning methods and algorithms, and have applied them to solve problems arising from diverse fields of sciences, engineering and humanities, ranging from genomics, neurosciences and health sciences to economics, finance and machine learning. For an overview, see \cite{BvdG2011} and \cite{HTW2015}.

Linear regression is often used to investigate the relationship between a response variable $Y$ and  explanatory variables $\bX = (X_1, \ldots, X_p)^\T$. In the high-dimensional linear model $Y=\bX^\T \bbeta^* + \varepsilon$, the coefficient $\bbeta^*$ is assumed to be sparse with support $S_0 = \supp(\bbeta^*)$. Variable selection techniques such as the forward stepwise regression, the Lasso \citep{TIB96} and folded concave penalized least squares \citep{FLi2001, ZL2008} are frequently used. However, it has been recently noted in \cite{FGH2012} that high dimensionality introduces large spurious correlations between response and unrelated covariates, which may lead to wrong statistical inference and false scientific discoveries. As an illustration, \cite{FSZ2015} considered a real data example using the gene expression data from the international `HapMap' project \citep{T2005}. There, the sample correlation between the observed and  post-Lasso fitted responses is as large as $0.92$. While conventionally it is a common belief that a correlation of 0.92 between the response and a fit is noteworthy, in high-dimensional scenarios, this intuition may no longer be true. In fact, even if the response and all the covariates are scientifically independent in the sense that $\bbeta^* = \mo$, simply by chance, some covariates will appear to be highly correlated with the response. As a result, the findings obtained via any variable selection techniques are hardly impressive unless they are proven to be better than by chance. To simplify terminology, in this paper we say that the discovery (by a variable selection method) is spurious if it is no better than by chance.

To guard against spurious discoveries, one naturally asks how good a response can be fitted by optimally selected subsets of covariates, even when the response variable and the covariates are not causally related to each other, that is, when they are independent. Such a measure of the goodness of spurious fit (GOSF) is a random variable whose distribution can provide a benchmark to gauge whether the discoveries by statistical machine learning methods any better than a spurious fit (chance).  Measuring such a goodness of spurious fit and estimating its theoretical distributions are the aims of this paper.  This problem arises from not only high-dimensional linear models and generalized linear models, but also robust regression and other statistical model fitting. To formally measure the degree of spurious fit, \cite{FSZ2015} derived the distributions of maximum spurious correlations, which provide a benchmark to assess the strength of the spurious associations (between response and independent covariates) and to judge whether discoveries by a certain variable selection technique are any better than by chance.

The response, however, is not always a quantitative value. Instead, it is often binary; for example, positive or negative, presence or absence and success or failure. In this regard, generalized linear models (GLIM) serve as a flexible parametric approach to modeling the relationship between explanatory and response variables \citep{MN1989}. Prototypical examples include linear, logistic and Poisson regression models which are frequently encountered in practice.

In GLIM, the relationship between the response and covariates is more complicated and cannot be effectively measured via Pearson correlation coefficient, which is essentially a measure of the linear correlation between two variables. We need to extend the concept of spurious correlation or the measure of goodness of spurious fit to more general models and study its null distribution. A natural measure of goodness of fit is the likelihood ratio statistic, denoted by ${\cal LR}_n(s,p)$, where $n$ is the sample size and $s$ is size of optimally fitted model.  It measures the goodness of spurious fit when $\bX$ and $Y$ are independent.  This generalization is consistent with the spurious correlation studied in \cite{FSZ2015}, that is, applying ${\cal LR}_n(s,p)$ to linear regression yields the maximum spurious correlation. We plan to study the limiting null distribution of $2 {\cal LR}_n(s,p)$ under various scenarios. This reference distribution then serves as a benchmark to determine whether the discoveries are spurious.

\begin{figure}[hbtp!]
  \centering
  \includegraphics[scale=0.4]{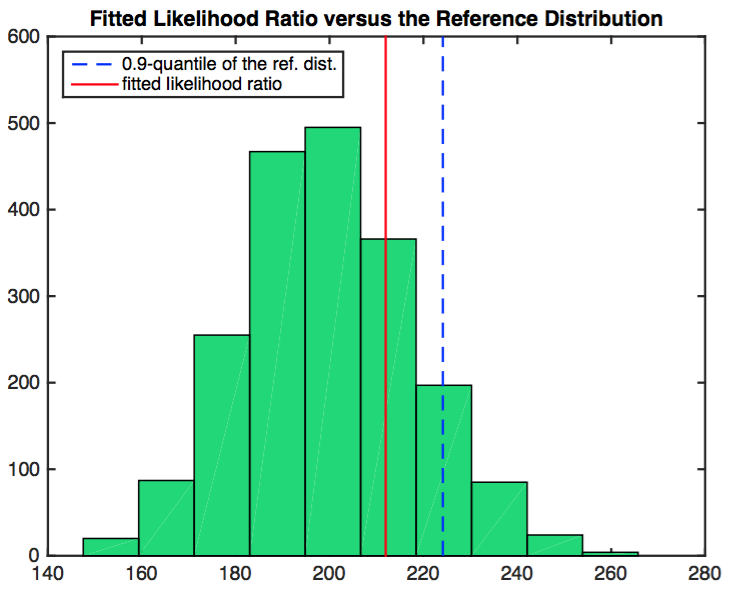}
  \caption{Lasso fitted likelihood ratio $2\hat{{\cal LR}}$ in comparison to
  the distribution of GOSF $2{\cal LR}_n(s,p)$ with $n=246$, $p=10,707$ and $s=40$. \label{fig0} }
\end{figure}

To gain further insights, let us illustrate the issue by using the gene expression profiles for $10,707$ genes from $251$ patients in the German Neuroblastoma Trials NB90-NB2004 \citep{O2006}. The response labeled as ``3-year event-free survival'' (3-year EFS) is a binary outcome indicating whether each patient survived 3 years after the diagnosis of neuroblastoma. Excluding five outlier arrays, there are $246$ subjects (101 females and 145 males) with 3-year EFS information available. Among them, 56 are positives and 190 are negatives. We apply Lasso using the logistic regression model with tuning parameter selected via ten-fold cross validation (40 genes are selected). The fitted likelihood ratio $2\hat{{\cal LR}} =  211.96$. To judge the credibility of the finding of these 40 genes, we should compare the value $211.96$ with the distribution of the Goodness Of Spurious Fit (GOSF) $2{\cal LR}_n(s,p)$ when $\bX$ and $Y$ are indeed independent, where $n=246$, $p=10,707$ and $s=40$. This requires some new methodology and technical work. Figure~\ref{fig0} shows the distribution of the GOSF estimated by our proposed method below and indicates how abnormal the value $211.96$ is.  It can be concluded that the goodness of fit to the binary outcome is not statistically significantly better than GOSF.

The above result shows that the 10-fold cross-validation chooses a too large model with 40 variables.  This prompts us to reduce the model sizes along the Lasso path such that their fits are better than GOSF.  The results are reported in Table~\ref{tab2}.
The largest model along the LASSO path that fits better than GOSF has model size 17.  We can use the cross-validation to select a model with model size no more than 17 or to select a best model among all models that fit better than GOSF.  This is another important application of our method.

\subsection{Structure of the paper}

In Section~\ref{sec2}, we introduce a general measure of spurious fit via generalized likelihood ratios, which extends the concept of spurious correlation in the linear model to more general models, including generalized linear models and robust linear regression. We also introduce a local adaptive majorization-minimization (LAMM) algorithm to compute the GOSF. Section~\ref{sec3} presents the main results on the limiting laws of goodness of spurious fit and their bootstrap approximations. For conducting inference, we use the proposed LAMM algorithm to compute the bootstrap statistic. In Section~\ref{sec4}, we discuss an application of our theoretical findings to high-dimensional statistical inference and model selection. Section~\ref{sec5} presents numerical studies. Proofs of the main results, Theorems~\ref{thm3.1} and \ref{thm3.3}, are provided in Section~\ref{sec6}; in each case, we break down the key steps in a series of lemmas with proofs deferred to the appendix.

\subsection{Notations}

We collect standard pieces of notation here for readers' convenience. For two sequences $\{a_n\}$ and $\{ b_n\}$ of positive numbers, we write
$a_n=O(b_n)$ or $a_n \lesssim b_n$ if there exists a constant $C>0$ such that
$ a_n /b_n \leq C $ for all sufficiently large $n$; we write $a_n \asymp b_n$
if there exist constants $C_1, C_2>0$ such that, for all $n$ large enough,
$C_1 \leq  a_n/ b_n \leq C_2$; and we write $a_n=o(b_n)$
if $\lim_{n\rightarrow \infty} a_n/b_n =0$,
respectively. For $a, b\in \bbr$, we write $a \vee b=\max(a, b)$.

For every positive integer $\ell$, we write $[\ell]= \{1, 2, \ldots, \ell\}$, and for any set $S$, we use $S^{{\rm c}}$ to denote its complement and $|S|$ for its cardinality. For any real-valued random variable $X$, its sub-Gaussian norm is defined by $\| X \|_{\psi_2} = \sup_{\ell \geq 1} \ell^{-1/2} (\e |X|^\ell)^{1/\ell}$. We say that a random variable $X$ is sub-Gaussian if $\| X \|_{\psi_2} <\infty$.

Let $p, q$ be two positive integers. For every $p$-vector $\bu = (u_1,\ldots, u_p)^\T$, we define its $\ell_q$-norm to be $\| \bu \|_q= \big(  \sum_{i =1}^p |u_i|^q  \big)^{1/q}$, and set $\| \bu \|_0=\sum_{i=1}^p I\{u_i \neq 0\}$. Let $\mathbb{S}^{p-1} = \{ \bu \in \bbr^p : \| \bu \|_2 =1  \}$ be the unit sphere in $\bbr^p$. Moreover, for each subset $S \subseteq [p]$ with $|S|=s \in [p]$, we denote by $\bu_S$ the $s$-variate sub-vector of $ \bu$ containing only the coordinates indexed by $S$. We use $\| \bM\|$ to denote the spectral norm of a matrix $\bM$.

\section{Goodness of spurious fit}
\label{sec2}

Let $Y, Y_1 , \ldots, Y_n $ be independent and identically distributed (i.i.d.) random variables with mean zero and variance $\sigma^2 > 0$, and $\bX, \bX_1, \ldots, \bX_n$ be i.i.d. $p$-dimensional random vectors. We write
$$
\bX = (X_1, \ldots, X_p)^\T , \, \bbX = (\bX_1, \ldots,  \bX_n)^\T \in \bbr^{n\times p}   \,\mbox{ and } \,   \bX_i = (X_{i1}, \ldots, X_{ip})^\T , \ \  i= 1,\ldots, n .
$$
For $s \in [p]$, the maximum $s$-multiple correlation between $Y$ and $\bX$ is given by
\begin{align}   \label{eq2.1}
	\hat{R}_n(s,p) = \max_{\balpha \in \bbr^p: \| \balpha \|_0 \leq s} \hat{\corr}_n(Y, \balpha^\T \bX ) ,
\end{align}
where $\hat{\corr}_n(\cdot, \cdot)$ denotes the sample Pearson correlation coefficient. When $Y$ and $\bX$ are independent, we regard $\hat{R}_n(s,p)$ as the maximum spurious (multiple) correlation. The limiting distribution of $\hat{R}_n(s,p)$ is studied in \cite{CJ2012} and \cite{FGH2012} when $s=1$ and $X \sim  N(\mo, \bI_p)$ (the standard normal distribution in $\bbr^p$), and later in \cite{FSZ2015} under a general setting where $s\geq 1$ and $\bX$ is sub-Gaussian with an arbitrary covariance matrix.

For binary data, the sample Pearson correlation is not effective for measuring the regression effect. We need a new metric. In classical regression analysis, the multiple correlation coefficient, also known as the $R^2$, is the proportion of variance explained by the regression model. For each submodel $S\subseteq [p]$, its $R^2$ statistic can be computed as
\begin{align}	\label{eq2.2}
	R^2_S
   =  \max_{\btheta \in \bbr^s}   \,\hat{\corr}_n^2  (Y, \bX_S^\T \btheta).
\end{align}
Then, the maximum $s$-multiple correlation $\hat{R}_n(s,p)$ can be expressed as the maximum $R^2$ statistic:
\begin{align} \label{eq2.3}
	\hat{R}_n^2(s,p) = \max_{S\subseteq [p] : |S| =s} R_S^2 .
\end{align}

The concept of $R^2$ can be extended to more general models. For binary response models, \cite{M1983} suggested the following generalization:
$
	-\log(1-R^2) = \frac{2}{n} \{ \ell(\hat{\bbeta}) - \ell(\mo) \},
$
where $\ell(\hat{\bbeta}) = \log L(\hat{\bbeta})$ and $\ell(\mo) = \log L(\mo)$ denote the log-likelihoods of the fitted and the null model, respectively. This motivates us to use the likelihood ratio as a generalization of the goodness of fit beyond the linear model.

Let $L_n(\bbeta)$, $\bbeta\in \bbr^p$ be the negative logarithm of a quasi-likelihood process of the sample $\{(Y_i,\bX_i )\}_{i=1}^n$. 
For a given model size $s\in [p]$, the best subset fit is $\hat{\bbeta}(s) := \argmin_{\bbeta \in  \bbr^p : \|\bbeta \|_0\leq s } L_n(\bbeta) $. The goodness of such a fit, in comparison with the baseline fit $L_n(\mo)$, can be measured by
\begin{align} \label{eq2.4}
 \mathcal{ LR}_n(s,p) := L_n(\mo) - L_n(\hat{\bbeta}(s)) =  L_n(\mo)  - \min_{\bbeta \in  \bbr^p : \|\bbeta \|_0\leq s } L_n(\bbeta).
\end{align}
When $\bX$ and $Y$ are independent, it becomes the Goodness OF Spurious Fit (GOSF).
According to \eqref{eq2.2} and \eqref{eq2.3}, this definition is consistent with the maximum spurious correlation when it is applied to the linear model with Gaussian quasi-likelihood, where $L_n(\bbeta ; \beta_0, \sigma) =  \frac{1}{2} \log( 2\pi \sigma^2 ) +  \frac{1}{2 } \|   \bY - \beta_0 - \bbX \bbeta \|_2^2  /\sigma^2$ and $\bY = ( Y_1, \ldots, Y_n)^\T$.

Throughout, we refer to $L_n(\cdot)$ as the loss function which is assumed to be convex. This setup encompasses the generalized linear models \citep{MN1989} with $L_n(\bbeta) = \sn \{ b(\bX_i^\T \bbeta) - Y_i \, \bX_i^\T \bbeta \}$ under the canonical link where $b(\cdot)$ is a model-dependent convex function (we take the dispersion parameter as one, as we don't consider the dispersion issue), robust regression with $L_n(\bbeta) = \sn |Y_i - \bX_i^\T \bbeta|$, the hinge loss $L_n(\bbeta) = \sn (1- Y_i \, \bX_i^\T \bbeta)_+$ in the support vector machine \citep{V1995} and exponential loss $L_n(\bbeta)  = \sn \exp(- Y_i \, \bX_i^\T \bbeta)$ in AdaBoost \citep{FS1997} in classification with $Y$ taking values $\pm 1$.

The prime goal of this paper is to derive the limiting laws of GOSF ${\cal LR}_n(s,p)$ in the null setting where the response $Y$ and the explanatory variables $\bX$ are independent. Here, both $s$ and $p$ can depend on $n$, as we shall use double-array asymptotics. We will mainly focus on the GLIM and robust linear regression that are of particular interest in statistics.


\subsection{Generalized linear models}
\label{GLIM.sec}

Recall that $(Y_1,\bX_1), \ldots, (Y_n, \bX_n)$ are i.i.d. copies of $(Y, \bX)$. Assume that the conditional distribution of $Y $ given $\bX =\bx \in  \bbr^p $ belongs to the canonical exponential family with the probability density function taking the form \citep{MN1989}
\begin{align}   \label{eq2.5}
	f(  y ; \bx  ,   \bbeta^*   ) =  \exp\big[  \{  y   \,\bx^\T \bbeta^* - b(\bx^\T \bbeta^*) \} / \phi + c(y, \phi) \big] ,
\end{align}
where $\bbeta^*=( \beta^*_1, \ldots, \beta^*_p)^\T$ is the unknown $p$-dimensional vector of regression coefficients, and $\phi>0$ is the dispersion parameter. 
The log-likelihood function with respect to the given data $\{(Y_i, \bX_i ) \}_{i=1}^n$ is $ \sn  c( Y_i, \phi)   + \phi^{-1}  \sn  \{ Y_i  \, \bX_i^\T \bbeta     - b(  \bX_i^\T \bbeta )  \}  $. For simplicity, we take $\phi=1$ with the exception that in the linear model with Gaussian noise, $\phi=\sigma^2$ is the variance. 
Two other showcases are
\begin{enumerate}

\item[1.]Logistic regression: $b(u) = \log(1+e^u)$, $u \in \bbr$ and $\phi=1$.

\item[2.]Poisson regression: $b(u) = e^u$, $u\in \bbr$ and $\phi=1$. 
\end{enumerate}

In GLIM, the loss function is  $L_n(\bbeta) = \sn \{ b(\bX_i^\T \bbeta ) -Y_i \, \bX_{i}^\T \bbeta \}$. By \eqref{eq2.4}, the generalized measure of goodness of fit for GLIM is
\begin{align}  \label{eq2.6}
	 \mathcal{ LR}_n(s,p) =  n b(0) - \min_{\bbeta \in \bbr^p: \| \bbeta \|_0 \leq s} L_n(\bbeta) .
\end{align}
In Section~\ref{sec3}, we derive under mild regularity conditions the limiting distribution of GOSF $ \mathcal{ LR}_n(s,p)$ in the null model. This extends the classical Wilks theorem \citep{W1938}. Here, we interpret $ \mathcal{LR}_n(s,p)$ as the degree of spuriousness caused by the high-dimensionality.

\subsection{$L_1$ regression}
\label{sec2.2}

In this section, we revisit the high-dimensional linear model
\begin{align}
	\bY = \bbX \bbeta^* + \beps   \ \ \mbox{ or } \ \  Y_i = \bX_i^\T \bbeta^* + \varepsilon_i, \, i=1,\ldots , n ,  \label{eq2.7}
\end{align}
where $\bY = (Y_1, \ldots, Y_n)^\T$ is the response vector and $\beps = (\varepsilon_1, \ldots, \varepsilon_n)^\T$ is the $n$-vector of measurement errors. Robustness considerations lead to least absolute deviation (LAD) regression and more generally quantile regression \citep{K2005}. For simplicity, we consider the $\ell_1$-loss $L_{n}(\bbeta) = \sn |Y_i - \bX_i^\T \bbeta|$, $\bbeta \in \bbr^p$. The generalized measure of goodness of fit \eqref{eq2.4} now becomes
\begin{align}  \label{eq2.8}
	 \mathcal{ LR}_n(s,p) =  \| \bY \|_1 - \min_{\bbeta \in \bbr^p: \| \bbeta \|_0 \leq s} L_n(\bbeta) .
\end{align}
The limiting distribution of GOSF $\mathcal{ LR}_n(s,p)$ is studied in Section~\ref{sec3.4}.

In particular, if $\varepsilon_1, \ldots, \varepsilon_n$ in \eqref{eq2.7} are i.i.d. from the double exponential distribution with the density $f_\varepsilon(u) = \frac{1}{2} e^{-|u| }$, $u\in \bbr$, the $\ell_1$-loss $L_{n}(\cdot)$ corresponds to the negative log-likelihood function. In general, we assume that the regression error $\varepsilon_i$ has median zero, that is, $\P(\varepsilon_i \leq 0)= \frac{1}{2}$. Hence, the conditional median of $Y_i$ given $\bX_i$ is $\bX_i^\T \bbeta^*$ for $i \in [n]$, and $\bbeta^* =  \argmin_{ \bbeta \in \bbr^p } \e_{\bX} \{ L_{n}(\bbeta) \}$, where $\e_{\bX}(\cdot) = \e(\cdot \, | \bX_1,\ldots, \bX_n)$ denotes the conditional expectation given $\{ \bX_i\}_{i=1}^n$.

\subsection{An LAMM algorithm}\label{sec2.3}

The computation of the best subset regression coefficient $\hbbeta(s)$ in \eqref{eq2.4} requires solving a combinatorial optimization problem with a cardinality constraint, and therefore is NP-hard. In the following, we suggest a fast and easily implementable method, which combines the forward selection (stepwise addition) algorithm and a local adaptive majorization-minimization (LAMM) algorithm \citep{LHY2000, FLSZ2015} to provide an approximate solution.

Our optimization problem is $\min_{ \bbeta \in \bbr^p: \| \bbeta \|_0\leq s} f(\bbeta)$, where $ f(\bbeta) = L_n(\bbeta)$. We say that a function $g(\bbeta \,|\, \bbeta^{(k)})$ majorizes $f(\bbeta)$ at the point $\bbeta^{(k)}$ if $f(\bbeta^{(k)})=g(\bbeta^{(k)} \,|\, \bbeta^{(k)})$ and $f(\bbeta) \leq g(\bbeta \,|\, \bbeta^{(k)})$ for all $\bbeta \in \bbr^p$. An majorization-minimization (MM) algorithm initializes at $\bbeta^{(0)}$ and then iteratively computes $\bbeta^{(k+1)} = \argmin_{\bbeta\in\bbr^p: \| \bbeta \|_0 \leq s} \, g(\bbeta \,|\, \bbeta^{(k)})$. The target value of such an algorithm is non-increasing since
\begin{equation}\label{eq2.9}
    f(\bbeta^{(k+1)}) \stackrel{\mbox{\tiny majorization}}{\leq} g(\bbeta^{(k+1)} \,|\, \bbeta^{(k)}) \stackrel{\mbox{\tiny minimization}}{\leq}
    g(\bbeta^{(k)}\,|\,\bbeta^{(k)}) \stackrel{\mbox{\tiny initialization}} {=} f(\bbeta^{(k)}).
\end{equation}

We now majorize $f(\bbeta)$ at $\hat{\bbeta}^{(k)}$ by an isotropic quadratic function
\begin{equation} \label{eq2.10}
	g_\lambda(\bbeta \, | \, \hat{\bbeta}^{(k)}) = f(\bbeta) + \Big\langle \nabla f(\hat \bbeta^{(k)}) , \bbeta - \hat{\bbeta}^{(k)} \Big\rangle + \frac{\lambda}{2} \| \bbeta - \hat{\bbeta}^{(k)} \|_2^2 , \ \ \bbeta \in \bbr^p .
\end{equation}
This is a valid majorization as long as $\lambda \geq \max_{\bbeta} \| \nabla^2 f( \bbeta)\|$ (this will be relaxed below).
The isotropic form on the right-hand side of (\ref{eq2.10}) allows a simple analytic solution given by
\begin{align}
 \hat{\bbeta}^{(k+1)}_\lambda =  \argmin_{\bbeta\in\bbr^p: \| \bbeta \|_0 \leq s} \, g(\bbeta \,|\, \bbeta^{(k)})  =  \big\{ \hat{\bbeta}^{(k)} -  \lambda^{-1}  \nabla f(\hat \bbeta^{(k)})   \big\}_{[1:s]}. \nn
\end{align}
Here, we used the notation that for any $\bbeta \in \bbr^p$, $\bbeta_{[1:s]}\in \bbr^p$ retains the $s$ largest (in magnitude) entries of $\bbeta$ and assigns the rest to zero.

\begin{remark}  \label{rmk1}
{\rm To implement the MM algorithm, we need to compute the gradient of the objective function of interest. In the $L_1$ regression, the loss function $L_{n}(\bbeta) = \sn |Y_i - \bX_i^\T \bbeta|$, $\bbeta \in \bbr^p$ is not  differentiable everywhere. Recall that the subdifferential of the absolute function $h(x) = |x|$, $x\in \bbr$ is given by
$$
	 \partial h(x)  = \left\{\begin{array}{ll}
	\{ 1 \},    & \mbox{if }x>0 ,  \\
	 $[$ -1, 1  $]$ ,   &  \mbox{if } x=0, \\
	\{ -1 \} , & \mbox{if } x <0.
	\end{array}  \right.
$$
With slight abuse of notation, we suggest a randomized algorithm using the stochastic subgradient
$
 \nabla L_{n}(\bbeta)  = \sn   I(Y_i - \bX_i^\T \bbeta>0) - I(Y_i - \bX_i^\T \bbeta>0) + U_i I(Y_i - \bX_i^\T \bbeta = 0) ,
$
where $U_1, \ldots, U_n$ are i.i.d.  random variables uniformly distributed on $[-1, 1]$.
}
\end{remark}

We propose to use the stepwise forward selection algorithm to compute an initial estimator $\hat{\bbeta}^{(0)}$.  As the MM algorithm decreases the target value as shown in (\ref{eq2.9}), the resulting target value is no larger than that produced by the stepwise forward selection algorithm.

To properly choose the isotropic parameter $\lambda>0$ without computing the maximum eigenvalue, we use the local adaptive procedure as in \cite{FLSZ2015}. Note that, in order to have a non-increasing target value, the majorization is not actually required. As long as  $f(\bbeta^{(k+1)}) \leq g(\bbeta^{(k+1)}\,|\,\bbeta^{(k)})$, arguments in (\ref{eq2.9}) hold. Starting from a prespecified value $\lambda = \lambda_0$, we successfully inflate $\lambda$ by a factor $\rho>1$. After the $\ell$th iteration, $\lambda = \lambda_\ell =  \rho^{\ell-1} \lambda_0$. We take the first $\ell$ such that $f( \hat \bbeta_{\lambda_\ell}^{(k+1)}  ) \leq g_{\lambda_\ell}( \hat \bbeta_{\lambda_\ell}^{(k+1) } \,|\, \hat \bbeta^{(k)})$ and set $\hat \bbeta^{(k+1) } = \hat{\bbeta}_{\lambda_\ell}^{(k+1)}$.  Such an $\ell$ always exists as a large $\ell$ will major the function $f$. We then continue with the iteration in the MM part. A simple criteria for stopping the iteration is that $|f(\hat \bbeta^{(k+1)} ) - f(\hat \bbeta^{(k)} )| \leq \epsilon$ for a sufficiently small $\epsilon$, say $10^{-5}$. We refer to \cite{FLSZ2015} for a detailed computational complexity analysis of the LAMM algorithm.

While the LAMM algorithm can be applied to compute $\hbbeta(s)$ in a general setting, in our application, the algorithm is mainly applied to compute GOSF under the null model (see Figure~\ref{fig0} and Section~\ref{sec3.5}).  From our simulation experiences, our algorithm delivers a good enough solution under the null model.  It always provides an upper certificate $f(\hbbeta_0)$ to the problem $\min_{ \| \bbeta \|_0\leq s} f(\bbeta)$, where $\hbbeta_0$ is the output of the LAMM algorithm.  As in \cite{BKM2016}, if needed to verify the accuracy of our method, a lower certificate is $f(\hbbeta_1)$, where $\hbbeta_1$ is the solution to the convex problem
$\min_{\| \bbeta \|_1 \leq B_s} f(\bbeta)$, and $B_s$ is a sufficient large constant so that the $L_0$-solution satisfies $\|\hbbeta(s)\|_1 \leq B_s$.  For example, under the null model, it is well known that $\|\hbbeta(s)\|_1 = O_{\P}\{ s \sqrt{(\log p)/n} \}$.  Therefore, we can take $B_s = C_s s \sqrt{(\log p)/n}$ for a sufficiently large constant $C_s$.  A data-driven heuristic approach is to take $ B_s = 2 \|\hbbeta_1(s)\|_1$ along the Lasso path such that  $\|\hbbeta_1(s)\|_0 = s$.

Note that the minimum target value falls in the interval $[f(\hbbeta_1), f(\hbbeta_0)]$.  If this interval is very tight, we have certified that $\hbbeta_0$ is an accurate solution.

\section{Asymptotic distribution of goodness of spurious fit}
\label{sec3}

\subsection{Preliminaries}
\label{sec3.1}

Define $p\times p$ covariance matrices
\begin{align} \label{eq3.1}
\bSigma = \e ( \bX \bX^\T)   \quad \mbox{ and } \quad    	\hat \bSigma= n^{-1} \sn \bX_i \bX_i^\T.
\end{align}
For $s \in [p]$, we say that $S \subseteq [p]$ is an $s$-subset if $|S|=s$. For every $s$-subset $S  \subseteq [p]$, let $\bSigma_{SS} $ and $\hat{\bSigma}_{SS}$ be the $s\times s$ sub-matrices of $\bSigma$ and $\hat{\bSigma}$ containing the entries indexed by $S \times S$, that is,
\begin{align}  	\label{eq3.2}
	   \bSigma_{S S} = \e ( \bX_{S}\bX_S^\T) ,  \quad  	\hat \bSigma_{SS } = n^{-1} \sn \bX_{iS} \bX_{iS}^\T .
\end{align}

\begin{assumption} \label{cond3.1}
{\rm The covariates are standardized to have unit second moment, that is, $\e(X_j^2)=1$ for $j=1,\ldots, p$. There exits a random vector $\bU  \in \bbr^p$ satisfying $\e( \bU \bU^\T  ) = \bI_p$, such that $\bX= \bSigma^{1/2} \bU$ and $A_0 := \sup_{\bv \in \mathbb{S}^{p-1}}\| \bv^\T \bU \|_{\psi_2} < \infty$. }
\end{assumption}

For $1\leq s\leq p$, the $s$-sparse condition number of $\bSigma$ is given by
\begin{align} \label{eq3.3}
		\gamma_s = \gamma_{s}(\bSigma) = \sqrt{ \lambda_{\max}(s)/ \lambda_{\min}(s)  }  ,
\end{align}
where $\lambda_{\max}(s)=\max_{\bu \in \mathbb{S}^{p-1} :  \| \bu \|_0 \leq s}   \bu^\T \bSigma \bu$ and $\lambda_{\min}(s)=\min_{\bu \in \mathbb{S}^{p-1} :  \| \bu \|_0 \leq s}  \bu^\T \bSigma \bu$ denote the $s$-sparse largest and smallest eigenvalues of $\bSigma$, respectively.

Let $\bG = ( G_1, \ldots, G_p)^\T \sim  N(\mo ,  \bSigma)$ be a centered Gaussian random vector with covariance matrix $\bSigma$. For any $s$-subset $S\subseteq [p]$, $\bG_S \sim N(\mo , \bSigma_{SS} )$. Define the random variable
\begin{align}
	R_0(s,p) = \max_{S\subseteq [p]: |S|=s} \|  \bSigma_{SS}^{-1/2} \bG_S \|_2 , \label{eq3.4}
\end{align}
which is the maximum of the $\ell_2$-norms of a sequence of dependent chi-squared random variables with $s$ degrees of freedom. The distribution of $R_0(s,p)$ depends on the unknown $\bSigma$ and can be estimated by the multiplier bootstrap in Section~\ref{sec3.5}.  It will be shown that this distribution is the asymptotic distribution of GOSF. In particular, for the isotropic case where $\bSigma = \bI_p$, $R_0(s, p) = G_{(1)}^2 + \cdots + G_{(s)}^2$, the sum of the largest $s$ order statistics of $p$ independent $\chi_1^2$ random variables.

\subsection{Generalized linear models}
\label{sec3.2}

For i.i.d. observations $\{(Y_i, \bX_i )\}_{i=1}^n$  from the distribution in \eqref{eq2.5}, define individual residuals $\varepsilon_i = Y_i - \e_{\bX}(Y_i) = Y_i -   b'(\bX_i^\T \bbeta^*)$ with conditional variance $\var_{\bX}(\varepsilon_i) = \phi b''(\bX_i^\T \bbeta^*)$, where $\var_{\bX}(\cdot)= \e_{\bX} \{ \cdot - \e_{\bX}(\cdot) \}^2$. In particular, under the null model, $Y$ is independent of $\bX$ with mean $\mu_Y  := \e(Y) =  b'(0)$ and variance $\sigma_Y^{ 2} := \var(Y )=\phi b''(0)$.

\begin{assumption} \label{cond3.3}
{\rm There exists $a_0>0$ such that $  \e \exp\{ u \sigma_Y^{-1} (Y-\mu_Y) \}   \leq  \exp( a_0 u^2/2 )$ holds for all $u	\in \bbr$. The function $b(\cdot)$ in \eqref{eq2.5} satisfies
\begin{align}
	 \min_{u: |u|\leq 1} b''(u) \geq a_1  \quad \mbox{ and } \quad   \max_{u : |u| \leq 1} | b'''(u) | \leq A_1    \label{eq3.5}
\end{align}
for some constants $a_1, A_1 >0 $. }
\end{assumption}

Condition~\ref{cond3.3} is satisfied by a wide class of GLIMs, including the logistic and Poisson regression models. The following theorem shows that, under certain moment and regularity conditions, the distribution of the generalized likelihood ratio statistic $2 {\cal LR}_n(s,p)$ can be consistently approximated by that of $ R_0^2(s,p)$ given in \eqref{eq3.4}.

\begin{theorem} \label{thm3.1}
Let Conditions~\ref{cond3.1} and \ref{cond3.3} be satisfied. Assume that $\phi=1$ in \eqref{eq2.5}, $p, n  \geq 3$ and $1\leq s \leq  \min(p,n)$. Then, under the null model \eqref{eq2.7} with $\bbeta^*=\mo$,
\begin{align}
	& \sup_{t\geq 0} \big| \P \big\{    2  {\cal LR}_n(s, p) \leq t  \big\}  - \P  \{   R_0^2(s,p)  \leq t \}  \big|  \nn \\
	& \qquad \qquad  \qquad \qquad \qquad \leq C \big[  \{ s  \log(\gamma_s pn) \}^{7/8}   n^{-1/8}   +   \gamma_s^{1/2}     \{ s \log(\gamma_s pn) \}^2  n^{-1/2}   \big] , \label{eq3.6}
\end{align}
where $C>0$ is a constant depending only on $a_0, a_1 , A_0, A_1$ in Conditions~\ref{cond3.1} and \ref{cond3.3}.
\end{theorem}

\begin{remark}
{\rm
We regard Theorem~\ref{thm3.1} as a nonasymptotic, high-dimensional version of the celebrated Wilks theorem. In the low-dimensional setting where $s=p$ is fixed, Theorem~\ref{thm3.1} reduces to the conventional Wilks theorem, which asserts that the generalized likelihood ratio statistic converges in distribution to $\chi_p^2$.  In addition, we also provide a Berry-Esseen bound in \eqref{eq3.6}.
}
\end{remark}

\subsection{Linear least squares regression}
\label{sec3.3}

As a specific case of GLIM, we consider the linear regression model \eqref{eq2.7} with the loss function $L_n(\bbeta) = \frac{1}{2} \| \bY - \bbX \bbeta\|_2^2$. The corresponding likelihood ratio statistic
\begin{align}  \label{eq3.7}
{\cal LR}_n(s,p)   =  \frac{1}{2} \| \bY \|_2^2  - \min_{\bbeta\in \bbr^p: \| \bbeta\|_0 \leq s} L_n(\bbeta)
\end{align}
then coincides with that in \eqref{eq2.6} with $b(u)=\frac{1}{2} u^2$.  We state the null limiting distribution of ${\cal LR}_n(s,p)$ in a general case, where $\varepsilon_1, \ldots, \varepsilon_n$ are i.i.d. copies of a sub-Gaussian random variable $\varepsilon$. Specifically, we assume that

\begin{assumption}  \label{cond3.2}
{\rm $\varepsilon$ is a centered, sub-Gaussian random variable with $\var(\varepsilon)=\sigma^2>0$ and $K_0:= \| \varepsilon \|_{\psi_2} < \infty$. Moreover, write $v_\ell = \e ( |\varepsilon|^\ell )$ for $\ell \geq 3$. }
\end{assumption}

The following corollary is a particular case of the general result Theorem~\ref{thm3.1} with $b(u)=\frac{1}{2} u^2$, $u\in \bbr$ and $\phi=\sigma^2$. By examining the proof of Theorem~\ref{thm3.1} and noting that $b'''  \equiv 0$, it can be easily shown that the second term on the right-side of \eqref{eq3.6} vanishes. Hence, the proof is omitted.

\begin{corollary}    \label{cor3.2}
Let Conditions~\ref{cond3.1} and \ref{cond3.2} hold. Assume that $p, n \geq 3$ and $1 \leq s \leq \min(p,n)$. Then, under the null model \eqref{eq2.7} with $\bbeta^* = \mo$,
\begin{align}
	  \sup_{t\geq 0} \big| \P\big\{  2  {\cal LR}_n(s,p)   \leq t \big\}  - \P \big\{  \sigma^2  R^2_0(s,p)  \leq t \big\}  \big|   \leq C   \{ s  \log(\gamma_s pn)  \}^{7/8}  n^{-1/8}    ,    \nn
\end{align}
where $C>0$ is a constant depending only on $A_0$ and $K_0$ in Conditions~\ref{cond3.1} and \ref{cond3.2}.
\end{corollary}

\begin{remark}  \label{rmk3.3}
{\em Under the null model, the variance $\sigma^2$ can be consistently estimated by $\hat{\sigma}_0^2 =  n^{-1} \sn (Y_i  -\bar{Y} )^2$, where $\bar{Y} = n^{-1} \sn Y_i$. Under the same conditions of Corollary~\ref{cor3.2}, it can be proved that
$$
	  \sup_{t\geq 0} \big| \P \big\{  2 {\cal LR}_n(s,p) \leq t \big\} -  \P\big\{ \hat{\sigma}_0^2 R^2_0(s,p)  \leq t \big\} \big|   \lesssim     \{ s  \log(\gamma_s pn)  \}^{7/8}  n^{-1/8} ,   \nn
$$
which is in line with Theorem~3.1 in \cite{FSZ2015}. To see this, note that
\begin{align}
2{\cal LR}_n(s,p)   & =   \| \bY \|_2^2  -  \min_{S\subseteq [p]: |S| =s} \min_{\btheta \in \bbr^s} \| \bY - \bbX_S \btheta \|_2^2 \nn \\
& =  \max_{S\subseteq [p]: |S| =s}   \bY^\T  \bbX_S (\bbX_S^\T \bbX_S)^{-1} \bbX_S^{\T} \bY =    \max_{\balpha \in \bbr^p: \| \balpha\|_0 \leq s }  (  \bY^\T \bbX \balpha ) ^2 / \|  \bbX  \balpha \|_2^2.  \nn
\end{align}
The estimator $\hat{\sigma}_0^2$, used in computing the maximum spurious correlation, can be seriously biased beyond the null model and hence adversely affect the power. Thus, we suggest using either the refitted cross-validation procedure \citep{FGH2012} or the scaled Lasso estimator \citep{SZ2012} to estimate $\sigma^2$.
}
\end{remark}

\subsection{Linear median regression}
\label{sec3.4}

We now state an analogous result to Theorem~\ref{thm3.1} regarding the $\ell_1$-loss considered in Section~\ref{sec2.2}.

\begin{assumption} \label{cond3.4}
{\rm The noise $\varepsilon_1,\ldots, \varepsilon_n$ in \eqref{eq2.7} are i.i.d. copies of a random variable $\varepsilon$ satisfying $\e  |\varepsilon |^\kappa  <\infty$ for some $1<\kappa \leq 2$. There exist positive constants $a_2<(\e |\varepsilon|)^{-1}$, $A_2$ and $A_3$ such that the distribution function $F_\varepsilon(\cdot)$ and the density function $f_\varepsilon(\cdot)$ of $\varepsilon$ satisfy
\begin{align}
	 2\max\{ 1- F_\varepsilon(u) ,  F_\varepsilon(-u)  \}  \leq   ( 1+ a_2 u )^{-1}  \ \ \mbox{ for all } u \geq 0,   \label{eq3.8} \\
	   \max_{u\in \bbr} f_\varepsilon(u)  \leq A_2  \quad \mbox{ and } \quad   \max_{u: |u|\leq 1} \max\big\{ | f'_\varepsilon(u+) | ,  |f'_\varepsilon(u-)|  \big\} \leq A_3.   \label{eq3.9}
\end{align}}
\end{assumption}

\begin{theorem} \label{thm3.3} If $p, n  \geq 3$ and $1\leq s \leq  \min(p,n)$, then under the null model \eqref{eq2.7} with $\bbeta^*=\mo$ and Conditions~\ref{cond3.1} and \ref{cond3.4}, we have
\begin{align}
	& \sup_{t\geq 0} \big| \P \big\{   2 {\cal LR}_n(s,p) \leq t  \big\}  - \P \big\{  R_0^2(s,p) / \{ 2 f_\varepsilon(0) \} \leq t \big\} \big|  \nn \\
	& \qquad \qquad  \qquad \qquad   \leq C_1 \, n^{1-\kappa } + C_2 \big[    \{ s  \log(\gamma_s pn)  \}^{7/8}  n^{-1/8}  +  \gamma_s^{1/4}     \{ s  \log(\gamma_s pn) \}^{3/2}  n^{-1/4}    \big],   \label{eq3.10}
\end{align}
where ${\cal LR}_n(s,p)$ is given by \eqref{eq2.8}, $C_1>0$ is a constant depending on $a_2$, $\kappa$, $\e  |\varepsilon|  $, $\e |\varepsilon|^\kappa $ and $C_2>0$ is a constant depending on $ a_2 , A_0, A_2$ and $A_3$ in Conditions~\ref{cond3.1} and \ref{cond3.4}.
\end{theorem}

\begin{remark}  \label{rmk3.2}
{\rm Under the null model, the unknown parameter $f_\varepsilon(0)$ can be consistently estimated by the kernel density estimator $\hat{f}_\varepsilon(0) = (nh)^{-1} \sn K( Y_i/h )$, where $K(\cdot)$ is a kernel function and $h=h_n>0$ is the bandwidth. For simplicity, we may use the Epanechnikov kernel function $K_{{\rm Epa}}(u)=\frac{3}{4}(1-u^2)I(|u|\leq 1)$ along with the rule-of-thumb bandwidth $h_{{\rm ROT}} = 2.34 \, \hat{\sigma}_0 n^{-1/5}$, where $\hat{\sigma}_0^2 = n^{-1} \sn (Y_i -\bar{Y})^2$. 
}
\end{remark}

\subsection{Multiplier bootstrap procedure}
\label{sec3.5}

The distribution of the random variable $R_0(s,p)$ given by \eqref{eq3.4} depends on the unknown covariance matrix $\bSigma$. In practice, it is natural to replace $\bSigma$ by $\hat{\bSigma}=n^{-1}\sn \bX_i \bX_i^\T$ and $\bG \sim  N(0, \bSigma)$ by $\hat{\bG} \sim  N(\mo ,  \hat{ \bSigma} )$ in the definition of $R_0(s,p)$.  With this substitution, the distribution of $R_0(s,p)$ can be simulated.  In particular, $\hat{\bG}$ can be simulated as $ n^{-1/2} \sn e_i \bX_i$, where $e_1,\ldots, e_n$ are i.i.d. standard normal random variables that are independent of $\{ \bX_i \}_{i=1}^n$. The resulting estimator is
\begin{align}
 R_n(s,p) = \max_{ S\subseteq [p]: |S| =s } \| \hat{\bSigma}^{-1/2}_{SS} \hat{\bG}_S \|_2 ,  \label{eq3.11}
\end{align}
which is a multiplier bootstrap version of $R_0(s,p)$. The following proposition follows directly from Theorem~3.2 in \cite{FSZ2015}.

\begin{proposition} \label{prop3.1}
Assume that Condition~\eqref{cond3.1} holds, $1\leq s \leq  \min(p,n)$ and $s\log(\gamma_s p n ) = o(n^{1/5})$ as $n\to \infty$. Then $\sup_{t\geq 0}  | \P \{ R_0(s,p) \leq t \} - \P \{   R_n(s,p)  \leq t\, | \bX_1, \ldots, \bX_n \}  | \rightarrow  0$ in probability.
\end{proposition}

The computation of $R_n(s,p)$ requires solving a combinatorial optimization. This can be alleviated by using the LAMM algorithm in Section~\ref{sec2.3}.  To begin with, by Remark~\ref{rmk3.3}, we write $R_n(s,p)$ in \eqref{eq3.11} as
\begin{equation}\label{eq3.12}
 R^2_n(s,p) =  \max_{S\subseteq [p]: |S|=s}  \mathbf{e}^\T \bbX_S   (\bbX_S^\T \bbX_S)^{-1} \bbX_S^\T \mathbf{e} =  \| \be \|_2^2 - \min_{ \bbeta   \in \bbr^p: \| \bbeta \|_0\leq s} \| \be - \bbX \bbeta \|_2^2,  \nn
\end{equation}
where $\be = (e_1, \ldots, e_n)^\T$ and $\bbX_S = (\bX_{1S}, \ldots, \bX_{nS})^\T$ for every subset $S\subseteq [p]$. 
This can be computed approximately by the LAMM algorithm in Section~\ref{sec2.3}, resulting in the solution $\hat{\bbeta}(s)$.
Finally, we set $R_n^2(s,p) = \| \be \|_2^2 - \| \be - \bbX \hat{\bbeta}(s) \|_2^2$.

The numerical performance may be improved by employing mixed integer optimization formulations \citep{BKM2016}. Such an attempt, however, is beyond the scope of the paper and we leave it for future research.

\section{Spurious discoveries and model selection}
\label{sec4}

Based on the theoretical developments in Section~\ref{sec3}, here we address the question whether discoveries by machine learning and data mining techniques for GLIM are any better than by chance. For simplicity, we focus on the Lasso. Let $q_\alpha(s, p)$ be the upper $\alpha$-quantile of the random variable $R_0(s,p)$ defined by \eqref{eq3.4}. Assume that the dispersion parameter $\phi$ in \eqref{eq2.5} equals 1. By Theorem~\ref{thm3.1}, we see that for any prespecified $\alpha \in (0, 1)$,
\begin{align}
	\P \big\{   2{\cal LR}_n(s,p)   \leq   q^2_\alpha(s,p)    \big\} \to 1- \alpha , \label{eq4.1}
\end{align}
where ${\cal LR}_n(s,p)$ is as in \eqref{eq2.6}.

Let $\hat{\bbeta}_\lambda = \argmin_{\bbeta} \{ L_n(\bbeta) + \lambda \| \bbeta \|_1 \} $ be the $\ell_1$-penalized maximum likelihood estimator with $ \hat{s}_\lambda = |\hat{S}_\lambda |=| \mathrm{supp} ( \hat{\bbeta}_\lambda )  |$, where $\lambda>0$ is the regularization parameter. The goodness of fit is likelihood ratio $L_n(\mo) - L_n(\hat{\bbeta}_\lambda )$. Since $\hat{s}_\lambda$ covariates are selected, it should be compared with the distribution of GOSF ${\cal LR}_n(s,p)$ by taking $s=\hat{s}_\lambda$. In view of \eqref{eq4.1}, if
\begin{align}
	 L_n(\hat{\bbeta}_\lambda )  \geq   L_n( \mo ) -  q_\alpha^2(\hat{s}_\lambda , p) / 2  =  n b(0) - q_\alpha^2(\hat{s}_\lambda , p)  /2 , \nn
\end{align}
then we may regard the discovery of variables $\hat{S}_\lambda$ as unimpressive, no better than fitting by chance, or simply spurious.

In practice,
the unknown quantile $q_\alpha( s, p)$ should be replaced by its bootstrap version $q_{n,\alpha}(s,p)$,  the upper $\alpha$-quantile of $R_n(s,p)$ defined by \eqref{eq3.11}. 
This leads to the following data-driven criteria for judging where the discovery $\hat{S}(\lambda)$ is spurious:
\begin{align}
	 L_n(\hat{\bbeta}_\lambda )  \geq    n b(0)  -  q^2_{n,\alpha}(\hat{s}_\lambda , p)  /2 .  \label{eq4.2}
\end{align}
The theoretical justification is given by Theorem~\ref{thm3.1} and Proposition~\ref{prop3.1}.  In particular, when the loss is quadratic, this reduces to the case studied by \cite{FSZ2015}.

The concept of GOSF and its theoretical quantile provide important guidelines for model selection.  Let $\hat{\bbeta}_{{\rm cv}}$ be a cross-validated Lasso estimator, which selects $\hat{s}_{{\rm cv}} = \| \hat{\bbeta}_{{\rm cv}} \|_0$ important variables. Due to the bias of the $\ell_1$ penalty, the Lasso typically selects far larger model size since the visible bias in Lasso forces the cross-validation procedure to choose a smaller value of $\lambda$. This phenomenon is documented in the simulations studies. See Table~\ref{tab1} in Section~\ref{sec5.2}. With an over-selected model, both the goodness of fit $\hat{{\cal LR}}_\lambda = L_n(\mo) - L_n(\hat \bbeta_\lambda)$ and the spurious fit can be very large, and so is the finite sample Wilks approximation error. To avoid over-selecting, we suggest an alternative procedure that uses the quantity $q_{n,\alpha}(s,p)$ as a guidance to choose the tuning parameter, which guards us from spurious discoveries. More specifically, for each $\lambda$ in the Lasso solution path, we compute $\hat{{\cal LR}}_\lambda $ and $q_{n,\alpha}(s,p) |_{s=\hat{s}_\lambda}$ with a prespecified $\alpha$. Starting from the largest $\lambda$, we stop the Lasso path the first time that the sign of $2\hat{{\cal LR}}_\lambda - q^2_{n,\alpha}( \hat{s}_\lambda ,p)$ is changed from positive to negative, and let $\hat{\lambda}_{{\rm fit}}$ be the smallest $\lambda$ satisfying $2\hat{{\cal LR}}_\lambda  \geq  q^2_{n,\alpha}( \hat{s}_\lambda ,p) $. Denote by $\hat{s}_{{\rm fit}}$ the corresponding selected model size. This value can be regarded as the maximum model size for Lasso (or any other variable selection technique such as SCAD) to choose from. Another viable alternative is to only select the best cross-validated model among those whose fit are better than GOSF. We will show in Section~\ref{sec5.2} by simulation studies that this procedure selects much smaller model size which is closer to the truth.

\section{Numerical studies}
\label{sec5}

\subsection{Accuracy of the Gaussian approximation}
\label{sec5.1}

First we ran a simulation study to examine how accurate the Gaussian approximation $R^2_0(s,p)$ is to the generalized likelihood ratio statistic $2  {\cal LR}_n(s,p)$ in the null model. To illustrate the method, we focus on the logistic regression model: $\P(Y=1|\bX) = \exp(\bX^\T \bbeta^*)/\{ 1+\exp(\bX^\T \bbeta^*) \}$.  Under the null model $\bbeta^* = 0$, $Y_1, \ldots, Y_n$ are i.i.d. Bernoulli random variables with success probability $1/2$.
Independent of $Y_i$'s, we generate $\bX_i \sim N(0, \bSigma)$ with two different covariance matrices: $\bSigma_1=(\rho^{|j-k|})_{1\leq j, k \leq p}$ and $\bSigma_2=( \sigma_{2, jk})_{1\leq j, k \leq p}$, where
$$
		\sigma_{2, jk } =    \big( \big|  |j-k| + 1 \big|^{2\rho}  +      \big|  |j-k|  - 1 \big|^{2\rho} - 2 |i-j|^{2\rho} \big)/2 , \ \ 1\leq j, k \leq p.
$$
The first design has an AR(1) correlation structure (a short-memory process), whereas the second design reflects strong long memory dependence. We take $\rho=0.8$ in both cases.

\begin{figure}[hbtp!]
  \centering
  \includegraphics[scale=0.4]{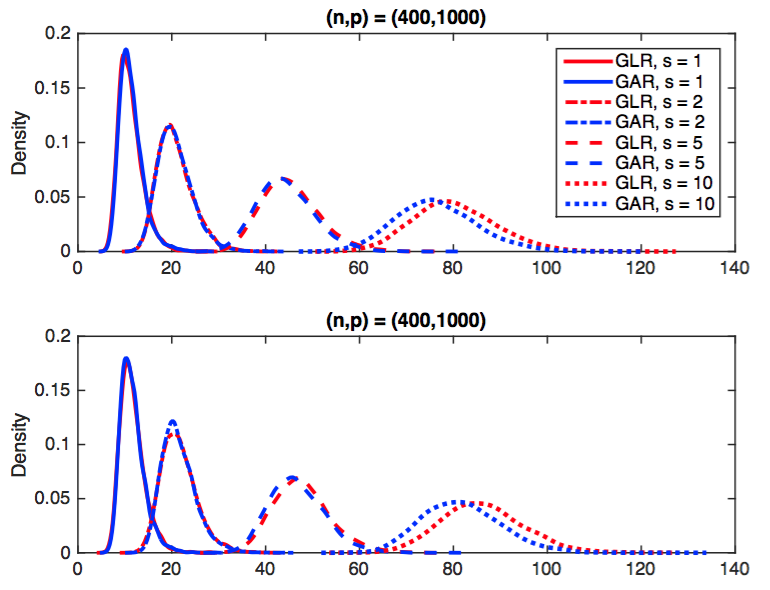}
  \caption{Distributions of generalized likelihood ratios (red) and Gaussian approximations (blue) based on 5000 simulations for $n=400$, $p = 1000$ and $s = 1, 2, 5, 10$ when $\bSigma$ is equal to $\bSigma_1$ (upper panel) or $\bSigma_2$ (lower panel).  \label{fig1} }
\end{figure}

Figure~\ref{fig1} reports the distributions of generalized likelihood ratios (GLRs) and their Gaussian approximations (GARs) when $n=400$, $p=1000$ and $s\in \{ 1, 2, 5, 10\}$. The results show that the accuracy of Gaussian approximation is fairly reasonable and is affected by the size of $s$ as well as the dependence between the coordinates of $\bX$.

\subsection{Detection of spurious discoveries}
\label{sec5.2}

In this section, we conduct a moderate scale simulation study to examine how effective the multiplier bootstrap quantile $q_{n,\alpha}(s,p)$ serves as a benchmark for judging whether the discovery is spurious. To illustrate the main idea, again we restrict our attention to the logistic regression model and the Lasso procedure.

The results reported here are based on 200 simulations with the ambient dimension $p = 400$ and the sample size $n$ taken values in $\{ 120,  160,  200 \}$. The true regression coefficient vector $\bbeta^* \in \bbr^p$ is $(3,-1,3,-1,3,0,\ldots, 0)^\T$. We consider two random designs: $\bSigma=\bI_p$ (independent) and $\bSigma = (0.5^{|j-k|})_{1\leq j,k\leq p}$ (dependent).

Let $\hat{\bbeta}_{{\rm cv}}$ be the five-fold cross-validated Lasso estimator, which selects a model of size $\hat{s}_{{\rm cv}} = \| \hat{\bbeta}_{{\rm cv}} \|_0$. For a given $\alpha\in (0,1)$, consider the spurious discovery probability (SDP)
\begin{align}
	\P\big\{   n \log(2)  - L_n(\hat{\bbeta}_{{\rm cv}} ) \leq   q^2_{n , \alpha}(\hat{s}_{{\rm cv}}  , p) /2\big\}, \nn
\end{align}
which is basically the probability of the type II error since the simulated model is not null. We take $\alpha = 0.1$ and compute the empirical SDP based on 200 simulations. For each simulated data set, $q_{n,\alpha}(s,p)|_{s=\hat{s}_{{\rm cv}} , \, p= 400}$ is computed based on 1000 bootstrap replications. The results are depicted in Table~\ref{tab1} below.

\begin{table}[h]
\centering
\caption{\label{size} The empirical power and the median size of the selected models with its robust standard deviation (RSD) in the parenthesis based on 200 simulations when  $p=400$ and $\alpha=10 \%$. RSD is the interquantile range divided by 1.34.}
\label{tab1}
\scriptsize{\begin{tabular}{ccccccccc}
\toprule
 & \multicolumn{2}{c}{$n=120$}& \multicolumn{2}{c}{$n=160$}  &   \multicolumn{2}{c}{$n=200$} \\
  \cline{2-3}\cline{4-5}\cline{6-7}   \vspace{-0.2cm}   \\
   &  Ind. & Dep. &  Ind. & Dep. &  Ind. & Dep.    \\ \midrule
 Power  & 0.595  & 0.750 & 0.925  & 0.980 & 1.000 & 1.000 \vspace{0.1cm}   \\
 $\hat{s}_{{\rm cv}}$ & 32.0 & 24.5 & 40.0 & 25.5 & 42.0 & 29.0     \\
 &    (13.43) & (11.94) & (13.81) & (12.69) & (14.18)  & (14.18)    \\
\bottomrule
\end{tabular}}
\end{table}\par

As reflected by Table~\ref{tab1}, the empirical power, which is one minus the empirical SDP, increases rapidly as the sample size $n$ grows. 
This is in line with our intuition that the more data we have, the less likely that the discovery by a variable selection method is spurious.
When the sample size is small, the SDP can be high and hence the discovery $\hat{S}_{{\rm cv}} = \supp(\hat{\bbeta}_{{\rm cv}})$ should be interpreted with caution. We need either more samples or more powerful variable selection methods.

We see from Table~\ref{tab1} that the Lasso with cross-validation selects far larger model size than the true one, which is 5. This is because the intrinsic bias in Lasso forces the cross-validation procedure to choose a smaller value of $\lambda$. 
We now use our procedure in Section~\ref{sec4} to choose the tuning parameter from the Lasso solution path. As before, we take $\alpha = 0.1$ in $q_{n,\alpha}(s, p)$ to provide an upper bound on the model size from perspective of guarding against spurious discoveries. The empirical median of $\hat{s}_{{\rm fit}}$ and its robust standard deviation are 9 and 1.87 over 200 simulations when $(n,p)=(200,400)$ and $\bSigma=  (0.5^{|j-k|})_{1\leq j,k\leq p}$. The feature over-selection phenomenon is considerably alleviated.

\subsection{Neuroblastoma data}
\label{sec5.3}

In this section, we apply the idea of detecting spurious discoveries to the neuroblastoma data reported in \cite{O2006}. This data set consists of 251 patients of the German Neuroblastoma Trials NB90-NB2004, diagnosed between 1989 and 2004. The complete data set, obtained via the MicroArray Quality Control phase-II (MAQC-II) project \citep{Shi2010}, includes gene expression over 10,707 probe sites. There are 246 subjects with 3-year event-free survival information available (56 positive and 190 negative). See \cite{O2006} for more details about the data sets.

For each $\lambda>0$, we apply Lasso using the logistic regression model to select $\hat{s}_\lambda$ genes. In particular, ten-fold cross-validated Lasso selects $\hat{s}_{{\rm cv}}  = 40$ genes. Then we calculate the goodness of fit
$ \hat{{\cal LR}}_\lambda :=  L_n(\mo) -  L_n(\hat{\bbeta}_\lambda)  = n \log (2) -  L_n(\hat{\bbeta}_\lambda)$. Along the Lasso path, we record in Table~\ref{tab2} the number of selected probes, the corresponding square-root the goodness of fit $( 2 \hat{{\cal LR}}_\lambda )^{1/2}  $ and upper $\alpha$-quantiles of the multiplier bootstrap approximations $R_0(s,p) |_{s=\hat{s}_\lambda, \, p = 10,707}$ with $\alpha = 10\%$ and $5\%$ based on 2000 bootstrap replications. For illustrative purposes, we only display partial Lasso solutions with selected model size $\hat{s}_\lambda$ lying between 20 and 40. From Table~\ref{tab2}, we observe that only the discovery of 17 probes has a generalized measure of the goodness of fit better than GOSF at $\alpha=5\%$, whereas the finding (of the 40 probes) via the cross-validation procedure is likely to over-select.

\begin{table}[h]
\centering
\caption{\label{size} Lasso fitted square-root likelihood ratio statistic, the mean cross-validated error, and upper $0.1$- and $0.05$-quantiles of the multiplier bootstrap approximation based on $2000$ bootstrap samples.}
\label{tab2}
\footnotesize{\begin{tabular}{ccccccc}
\toprule
  $\lambda$ & $\hat{s}_\lambda$  &  $ (2  \hat{{\cal LR}}_\lambda  )^{1/2}  $    &  $q_{n, 0.1}(\hat{s}_\lambda, p)  $   & $ q_{n, 0.05}(\hat{s}_\lambda, p)  $  & Mean Cross-Validated Error \\ \midrule
$ 0.2117 $ & 3 & 9.1389  & 6.4898 &  6.6519 &  1.0641 \\ \midrule
$ 0.1929 $ & 4 & 9.4753  & 7.2464 &  7.4353 & 1.0450 \\ \midrule
$ 0.1841 $ & 6 & 9.7273  & 8.4241 &  8.6061  & 1.0346 \\ \midrule
$  0.1678 $ & 7 & 10.1670 & 8.8959 &  9.0750 & 1.0092 \\ \midrule
$ 0.1601 $ & 8 & 10.3675  & 9.3121 &  9.5102 & 0.9974 \\ \midrule
$ 0.1459 $ & 9 & 10.7263  & 9.7115 &  9.9097 & 0.9751 \\ \midrule
$ 0.1329  $ & 11 & 11.0739 & 10.3954  & 10.6071 &  0.9543 \\ \midrule
$ 0.1269    $ & 12 & 11.2376  & 10.7042 &  10.9207 & 0.9452 \\ \midrule
$ 0.1211    $ & 13 & 11.4330  & 10.9875 &  11.2085 & 0.9359  \\ \midrule
$ 0.1104    $ & 14 & 11.7764  & 11.2576 &  11.4849 & 0.9186 \\ \midrule
$ 0.1006    $ & 15 & 12.0756 & 11.5084 &  11.7407 & 0.9006 \\ \midrule
$ 0.0960   $ &  17 & \textbf{12.2096} & \textbf{11.9664} & \textbf{12.2000} &  0.8934 \\ \midrule
$  0.0875 $ & 20 & \textbf{12.4788} & \textbf{12.5543} &  \textbf{12.7891} & 0.8815 \\ \midrule
$ 0.0761  $ &  25 & 12.9535 & 13.3824 & 13.6022 & 0.8651 \\ \midrule
$  0.0575  $&  31 & 13.8675 & 14.1407 &  14.3703 & 0.8361 \\ \midrule
$ 0.0456  $ & 40 &  14.5588 &  14.9712 & 15.2099  & 0.8255 \\
\bottomrule
\end{tabular}}
\end{table} \par

\section{Proofs}
\label{sec6}

We now turn to the proofs of Theorems~\ref{thm3.1} and \ref{thm3.3}. In each proof, we provide the primary steps, with more technical details stated as lemmas and proved in the appendix.


\subsection{Proof of Theorem~\ref{thm3.1}}

Throughout, we work with the quasi-likelihood ${\cal L}_n(\bbeta) = - L_n(\bbeta) = \sn  \{ Y_i \, \bX_i^\T \bbeta - b(\bX_i^\T \bbeta)  \}$ and consider the general case where the dispersion parameter $\phi$ in \eqref{eq2.5} is specified (not necessarily equals 1 to facilitate the derivations for the normal case). For a given $s\in [p]$, define
\begin{align}
   Q_{n}(s, p)  = \max_{ \bbeta \in  \bbr^p : \| \bbeta \|_0 \leq s   } {\cal L}_n( \bbeta)  \quad \mbox{ and } \quad Q_n^* = {\cal L}_n( \mo ).  \nn
\end{align}
We divide the proof into three steps. First, for each $s$-subset $S\subseteq [p]$, we prove Wilks's result for the $S$-restricted model where only a subset of the covariates indexed by $S$ are included. Specifically, we show that the square root deviation of the $S$-restricted maximum log-likelihood from its baseline value under the null model can be well approximated by the $\ell_2$-norm of the normalized score vector. Second, based on a high-dimensional invariance principle, we prove the Gaussian/chi-squared approximation for the maximum of the $\ell_2$-norms of normalized score vectors. Finally, we apply an anti-concentration argument to construct non-asymptotic Wilks approximation for $2\{  Q_{n}(s, p) - Q_n^* \}$.

\medskip
\no
{\it \textbf{Step~1: Wilks approximation}}. In the null model where $Y$ and $\bX$ are independent, the true parameter $\bbeta^*$ in \eqref{eq2.5} is zero, and thus the density function of $Y$ has the form $f(y )=\exp\{ - \phi^{-1}   b(0) + c(y,\phi) \}$. Moreover, we have
$$
\argmax_{  \bbeta  \in  \bbr^p }  \e_{\bX}  \{ {\cal L}_n( \bbeta ) \} =\argmax_{ \bbeta \in  \bbr^p } \sn \e_{\bX} \{ Y_i \,\bX_i^\T \bbeta -  b(\bX_i^\T \bbeta) \} = \mo .
$$
To this see, note that in model \eqref{eq2.5} with $\bbeta^*= \mo $, $\e( Y )= b'(0)$ and $\var(Y) =\phi b''(0)$.  This implies that $\e_{\bX} \{ \mathcal{L}_n( \bbeta ) \}  =  \sn \{  b'(0) \bX_i^\T \bbeta - b(\bX_i^\T \bbeta) \}$.
This function is strictly concave with respect to $\bbeta$ and $\bbeta = \mo$ satisfies its first order condition, and hence is its maximizer.

For each $s$-subset $S \subseteq [p]$, define the $S$-restricted log-likelihood ${\cal L}_n^S(   \btheta  )    = \sn \{ Y_i \, \bX_{iS}^\T \btheta - b(\bX_{iS}^\T \btheta) \}$ and the score function $\nabla {\cal L}_n^S( \btheta)  = \sn  \{ Y_i  - b'( \bX_{i S}^\T \btheta) \}\bX_{i S}$, $\btheta \in \bbr^s$. In this notation, it can be seen from \eqref{eq2.6} that
\begin{align}
	Q_{n}(s, p) = \max_{S\subseteq[p] : |S|=s} \max_{  \btheta\in  \bbr^s } {\cal L}_n^S( \btheta ) =  \max_{S\subseteq[p] : |S|=s} {\cal L}_n^S( \hat{\btheta}_S )  ,  \label{Qn.equiv}
\end{align}
where
\begin{align}
	\hat{\btheta}_S  = (  \hat{\theta}_{S,1}  , \ldots , \hat{\theta}_{S,s}  )^\T = \argmax_{ \btheta \in  \bbr^s } {\cal L}_n^S( \btheta )  \label{mle.sub-model}
\end{align}
denotes the maximum likelihood estimate of the target parameter for the $S$-restricted model, which is given by $\btheta_S^* : = \argmax_{ \btheta \in \bbr^s }  \e_{\bX}  \{ {\cal L}^S_n( \btheta) \} =  \mo$.

Given the i.i.d. observations $\{(Y_i,\bX_i )\}_{i=1}^n$, $\nabla \e_{\bX} \{ {\cal L}^S_n( {\btheta}) \} = \sn \{  b'(0)  -    b'( {\bX}_{i S}^\T {\btheta} ) \} {\bX}_{i S} $ and $\bH_S(\btheta) := - \nabla^2 \e_{\bX} \{ {\cal L}^S_n(\btheta) \}  =\sn b''(\bX_{i S}^\T \btheta) \bX_{i S} \bX_{i S}^\T$ for $\btheta\in \bbr^s$. In particular, write
\begin{align}
	 \bH^*_{S} : = \bH_S( \mo )  =   n  b''(0) \, \hat{\bSigma}_{SS}   \label{infor.mat.sub-model}
\end{align}
for $\bSigma_{SS}$ as in \eqref{eq3.2}. Further, define the $S$-restricted normalized score
\begin{align}
	\hat \bxi_S  =  \bH_S^{* -1/2} \nabla   {\cal L}^S_n(  \mo )  =  \{n b''(0) \}^{-1/2} \, \hat{\bSigma}_{SS}^{-1/2} \sn  \varepsilon_i \bX_{i S}, \qquad
 \varepsilon_i = Y_i - b'(0). \label{normalized.score.sub-model}
\end{align}

The following result is a conditional analogue of Corollary~1.12 in the supplement of \cite{S2012}, which provides an exponential inequality for the $\ell_2$-norm of $\hat \bxi_S$ given $\{ \bX_i \}_{i=1}^n$. The proofs of this Lemma and other lemmas can be found in the appendix. 

\begin{lemma} \label{quadratic.concentration.lemma}
Assume that Conditions~\ref{cond3.1} and \ref{cond3.3} hold. Then, for every $t \geq 0$,
\begin{align}
	\P_{\bX} \big\{  \| \hat \bxi_S \|_2^2  \geq a_0 \phi \Delta(s, t )  \big\}  \leq 2 e^{-t}  \label{quadratic.concentration.ineq}
\end{align}
holds almost surely on the event $\{\hat{\bSigma}_{SS} \succ \mo \}$, where
\begin{align}  \label{Delta.st}
 \Delta(s, t ) := \begin{cases}
  s +  (8ts)^{1/2} , \quad & {\rm if } ~0\leq t \leq \frac{1}{18} (2s)^{1/2}, \\ s +6t  , \quad & {\rm if } ~ t > \frac{1}{18}  (2s)^{1/2}.
\end{cases}
\end{align}
\end{lemma}

The following lemma characterizes the Wilks phenomenon from a non-asymptotic perspective. Recall that $\hat{\btheta}_S$ at \eqref{mle.sub-model} is the $S$-restricted maximum likelihood estimator, and in the null model, ${\cal L}^S_n( \mo  ) = {\cal L}_n( \mo  ) =  -n   b(0)$, $\sigma_Y^{2} = \var(Y) = \phi b''(0)$. For every $\tau >0$, define the event
\begin{align}
	\mathcal{E}_0( \tau ) =   \bigcap_{S\subseteq [p]: |S|=s}  \bigg\{ \hat{\bSigma}_{SS}  \succ \mo ,  \,\max_{1\leq i\leq n }  \bX_{i S}^\T \hat{\bSigma}_{SS}^{-1} \bX_{i S}  \leq  \tau  \bigg\} . \label{event1.def}
\end{align}

\begin{lemma} \label{Wilks.approximation}
Assume that Conditions~\ref{cond3.1} and \ref{cond3.3} hold. Then,  on the event $\mathcal{E}_0(\tau)$, for any $\tau>0$,
\begin{align}
 \P_{\bX} \bigg( \max_{S\subseteq [p] : |S| = s}  \Big| [ 2 \{ {\cal L}_n^S(\hat{\btheta}_S) - {\cal L}_n(\mo) \}  ]^{1/2}  -   \|	\hat \bxi_S \|_2   \Big|  \leq C_1 \, \phi \tau^{1/2} \frac{s\log(pn)}{\sqrt{n}}   \, \bigg) \leq 5 n^{-1}  	\label{Wilks.concentration}
\end{align}
whenever $n \geq C_2  \,   \phi \tau s\log(pn)$, where $C_1$ and $C_2$ are positive constants depending only on $a_0, a_1, A_1$ and $b''(0)$.

\end{lemma}

To apply Lemma~\ref{Wilks.approximation}, we need to show first that for properly chosen $\tau $, the event $ \mathcal{E}_0( \tau )$ occurs with high probability. First, applying Theorem 5.39 in \cite{V2012} to the random vectors $\bSigma_{SS}^{-1/2} \bX_{1S}, \ldots, \bSigma^{-1/2}_{SS} \bX_{nS}$ yields that, for every $t\geq 0$,
\begin{align}
	\big\| \bSigma_{SS}^{-1/2} \hat{\bSigma}_{SS}  \bSigma_{SS}^{-1/2} - \bI_s \big\| = \big\| n^{-1}\bSigma_{SS}^{-1/2}\bbX_S^\T \bbX_S \bSigma_{SS}^{-1/2} - \bI_s \big\| \leq \max(\delta, \delta^2)  \label{hatS.positivity}
\end{align}
holds with probability at least $1-2e^{-t}$, where $\delta = C_3 (s \vee t)^{1/2}n^{-1/2}  $, and $C_3>0$ is a constant depending only on $A_0$. This, together with Boole's inequality implies by taking $t=s\log\frac{ep}{s} + \log n$ that, with probability at least $1-2n^{-1}$,
\begin{align}
	\max_{S\subseteq [p]: |S|=s} \big\| \bSigma_{SS}^{-1/2} \hat \bSigma_{SS} \bSigma_{SS}^{-1/2} - \bI_{s } \big\| \leq C_3 \bigg( \frac{s\log\frac{ep}{s} + \log n}{n} \bigg)^{1/2} \leq \frac{1}{2} \label{uniform.hatS.positivity}
\end{align}
whenever $n \geq 4 C_3^2 (  s\log \frac{ep}{s} + \log n )$. Providing \eqref{uniform.hatS.positivity} holds, the smallest eigenvalue of $\bSigma_{SS}^{-1/2} \hat \bSigma_{SS} \bSigma_{SS}^{-1/2}$ is bounded from below by $\frac{1}{2}$ so that $\lambda_{\min}(\hat \bSigma_{SS} )  \geq \frac{1}{2} \lambda_{\min}(\bSigma_{SS})$. Moreover,
\begin{align}
	& \bX_{iS}^\T \hat \bSigma_{SS}^{-1} \bX_{iS}
\leq 2 \lambda^{-1}_{\min}(\bSigma_{SS})   \| \bX_{iS} \|_2^2 \leq 2s  \lambda^{-1}_{\min}(\bSigma_{SS})    \max_{j \in S} X_{ij}^2 .
\label{quadratic.ubd1}
\end{align}
For the last term on the right-hand side of \eqref{quadratic.ubd1}, let $\be_j = (0, \ldots, 0, 1, 0, \ldots, 0)^\T$ be the unit vector in $\bbr^p$ with 1 at the $j$th position and note that $X_{ij} = \be_j^\T \bX_i = \be_j^\T \bSigma_{SS}^{1/2} \bU_i  $ with $\| \be_j^\T \bSigma^{1/2} \|_2=1 $, where $\bU_1,\ldots, \bU_n$ are i.i.d. $p$-dimensional random vectors with covariance matrix $\bI_p$. By Condition~\ref{cond3.1}, $\|  X_{ij}   \|_{\psi_2}   =   \| \be_j^\T \bSigma^{1/2} \bU_i  \|_{\psi_2} \leq  A_0 $ and hence for every $t\geq 0$,
\begin{align}
	 \P\bigg(  \max_{1\leq i\leq n}\max_{1\leq j\leq p} X_{ij}^2 \geq t \bigg)  \leq   2\sn \sump \exp(  - C_4^{-1} t  ) \leq 2 \exp\{\log(pn) - C_4^{-1} t \}, \nn
\end{align}
where $C_4>0$ is a constant depending only on $A_0$. This, together with \eqref{quadratic.ubd1} implies by taking $t=2 C_4 \log(pn)$ that, with probability at least $1- 3n^{-1}$,
\begin{align}
\max_{1\leq i\leq n}\max_{S\subseteq [p]: |S|=s} \bX_{iS}^\T \hat \bSigma_{SS}^{-1} \bX_{iS}  \leq 2 \lambda^{-1}_{\min}(s )    \{ 1 + 2C_4 \, s\log(pn)  \}. \label{quadratic.ubd2}
\end{align}

Now, by \eqref{event1.def} and \eqref{quadratic.ubd2}, we take $\tau_0 = 2 \lambda_{\min}^{-1}(s )  \{ 1 + 2C_4  \, s\log(pn)  \}$ such that the event $ \mathcal{E}_0(\tau_0)$ occurs with probability greater than $1-3 n^{-1}$  as long as $n \geq 4 C_3^2  (  s\log \frac{ep}{s} + \log n )$. This, together with Lemma~\ref{Wilks.approximation} yields that with probability at least $1- 8 n^{-1}$,
\begin{align}
	\max_{S\subseteq[p]: |S|=s} \Big|  [ 2  \{ {\cal L}_n^S(\hat{\btheta}_S) - {\cal L}_n(\mo ) \}   ]^{1/2}  -  \|	\hat \bxi_S \|_2 \Big|  \leq C_5 \,  \phi \lambda^{-1/2}_{\min}(s)     \{s\log(pn)\}^{3/2}  n^{-1/2}   \label{uniform.Wilks.approxi}
\end{align}
whenever $n \geq C_6   (1\vee \phi)  \lambda^{-1}_{\min}(s) \{s\log(pn)\}^2$, where $C_5, C_6 >0$ are constants depending only on $a_0,  a_1 , A_0,  A_1$ and $b''(0)$.

\medskip
\no
{\it \textbf{Step~2: Gaussian approximation}}. For any $i=1,\ldots, n$ and $S \subseteq [p]$, define $\bZ_i = \{b''(0)\}^{-1/2} \varepsilon_i \bX_i $ and $\bZ_{iS} =  \{b''(0)\}^{-1/2} \varepsilon_i \bX_{iS}$ such that $ \hat{\bxi}_S = n^{-1/2} \sn \hat \bSigma_{SS}^{-1/2} \bZ_{iS} $. Moreover, define
\begin{align}
	\bxi = n^{-1/2} \sn \bZ_i \quad  \mbox{ and }  \quad  \bxi_S = n^{-1/2} \sn  \bSigma_{SS}^{-1/2} \bZ_{iS}.  \label{xiS.def}
\end{align}
The following result shows that for each $s$-subset $S\subseteq [p]$, the $\ell_2$-norm of the $S$-restricted normalized score $   \hat{\bxi}_S$ is close to that of $ \bxi_S$ with overwhelmingly high probability.

\begin{lemma} \label{GAR.lem1}
Assume that Condition~\ref{cond3.1} holds. Then, for every $s$-subset $S\subseteq [p]$ and for every $0\leq t\leq \frac{3}{4} ( n - 2s )$,
\begin{align}
	\P \Big[ \Big|   \| \hat{\bxi}_S \|_2  - \|  \bxi_S \|_2  \Big| > C_7  \{ (s +t) \phi \Delta(s,t) \}^{1/2} n^{-1/2}  \Big] \leq 12.4\,e^{-t}, \label{quadratic.approxi}
\end{align}
provided that $n\geq C_8 (s +t) $, where $\Delta(s,t)$ is as in \eqref{Delta.st} and $C_7, C_8 >0$ are constants depending only on $a_0$ and $A_0$.
\end{lemma}

Using the union bound and taking $t=s\log \frac{ep}{s} + \log n$ in Lemma~\ref{GAR.lem1}, we see that with probability at least $1-12.4\, n^{-1}$,
\begin{align}
\max_{S\subseteq[p]: |S|=s} \big|  \| \hat{\bxi}_S \|_2  - \|  \bxi_S \|_2  \big|  \leq C_7 \,  \phi^{1/2}  (  s\log \tfrac{ep}{s} + \log n )  \, n^{-1/2}  \label{unif.quadratic.approxi}
\end{align}
whenever $n\geq C_9  ( s \log \frac{ep}{s} + \log n )$.

Note that, the random vectors $\bxi$ and $\bxi_S, S\subseteq[p]$ defined in \eqref{xiS.def} satisfy $\e (\bxi) = \mo $, $\e( \bxi \bxi^\T) = \phi \bSigma$, $\e(\bxi_S) = \mo $ and $\e(\bxi_S \bxi_S^\T) = \phi \bI_s$. The following lemma provides a coupling inequality, showing that the random variable $\max_{S\subseteq [p]: |S|=s}  \| \phi^{-1/2}\bxi_S \|_2$ can be well approximated, with high probability, by some random variable which is distributed as the maximum of the $\ell_2$-norms of a sequence of normalized Gaussian random vectors, that is, $\{ \|  \bSigma_{SS}^{-1/2} \bG_S \|_2 : S\subseteq [p], |S| =s \}$.

\begin{lemma} \label{GAR.lem2}
Assume that Condition~\ref{cond3.1} holds. Then, there exists a random variable $T_0 \stackrel{d}{=}  R_0(s,p) $ such that for any $\delta\in (0, 1]$,
\begin{align}
	\bigg| \max_{S\subseteq [p]: |S| =s }  \| \phi^{-1/2} \bxi_S \|_2 - T_0 \bigg| \leq C_{10}  \big[ \delta +  \{ s\log(\gamma_s pn ) \}^{1/2}  n^{-1/2} +  \{s\log(\gamma_s pn )\}^2 n^{- 3/2}  \big]				\label{GAR.BE}
\end{align}
holds with probability greater than $1-C_{11}  \big[   \delta^{-3}n^{-1/2} \{s \log(\gamma_s pn)\}^2 \, \vee \, \delta^{-4}n^{-1} \{s \log(\gamma_s pn)\}^5 \big] $, where $C_{10}, C_{11}>0$ are constants depending only on $a_0$ and $A_0$
\end{lemma}

\noindent
{\it \textbf{Step~3: Completion of the proof}}. We now apply an anti-concentration argument to construct the Berry-Esseen bound for the square root of the excess $ 2\phi^{-1} \{ Q_{n}(s, p) - Q_n^*\}$. To this end, taking $\delta=\{s\log(\gamma_s pn)\}^{3/8}  n^{-1/8} $ in Lemma~\ref{GAR.lem2} leads to that, with probability at least $1-C_{11}   \{s\log(\gamma_s pn)\}^{7/8}  n^{-1/8} $,
\begin{align}
\bigg| \max_{S\subseteq [p]: |S| =s }  \| \phi^{-1/2} \bxi_S \|_2 - T_0 \bigg|  \leq C_{12}  \{ s\log(\gamma_s pn) \}^{3/8} n^{-1/8}  \label{uniform.coupling}
\end{align}
whenever $n \geq  \{s\log(\gamma_s pn)\}^3$. Further, for $R_0(s,p)$ in \eqref{eq3.4}, note that
\begin{align}
R^2_0(s,p) = \max_{S\subseteq[p]: |S| =s} \max_{\bu \in \mathbb{S}^{s-1}} \frac{ ( \bu^\T \bG_S )^2 }{ \bu^\T \bSigma_{SS} \bu }	 = \max_{\bu \in \mathcal{F}(s,p)} \frac{ ( \bu^\T \bG)^2 }{ \bu^\T \bSigma  \bu }, \nn
\end{align}
where $\bG \sim N(\mo , \bSigma )$ and $\mathcal{F}(s,p) :=   \{ \bx \mapsto  \bu^\T \bx   : \bu  \in \mathbb{S}^{p-1} ,  \|  \bu \|_0 \leq s  \}$ is a class of linear functions $\bbr^p \mapsto \bbr$. Hence, it follows from Lemma~7.3 in \cite{FSZ2015} with slight modification and Lemma~A.1 in the supplement of \cite{CCK2014} that, for every $t>0$,
\begin{align}
	\sup_{u\geq 0} \P ( | T_0  - u | \leq t  ) = \sup_{u\geq 0} \P \big\{ | R_0(s,p) - u | \leq t  \big\} \leq C_{13} \big( s\log \tfrac{\gamma_s e p}{s} \big)^{1/2} t , \label{R*.anti-concentration}
\end{align}
where $C_{13}>0$ is an absolute constant. Combining \eqref{R*.anti-concentration} with the preceding results \eqref{uniform.Wilks.approxi}, \eqref{unif.quadratic.approxi} and \eqref{uniform.coupling} proves \eqref{eq3.6}.  \qed

\subsection{Proof of Theorem~\ref{thm3.3}}

The main strategy of the proof is similar to that of Theorem~\ref{thm3.1} but technical details are substantially different. As before, we define the quasi-likelihood ${\cal L}_n(\bbeta) = - \sn |Y_i - \bX_i^\T \bbeta|$, $\bbeta \in \bbr^p$, and observe that
$ \max_{\bbeta \in  \bbr^p: \| \bbeta \|_0\leq s} {\cal L}_n(\bbeta)  = \max_{S\subseteq [p]: |S|=s} \max_{\btheta \in \bbr^s} {\cal L}_n^S(\btheta)$, where ${\cal L}_n^S(\btheta) = - \sn |Y_i - \bX_{iS}^\T \btheta |$. In the null model \eqref{eq2.7} with $\bbeta^* = \mo $, we have for each $s$-subset $S \subseteq [p]$, $ \argmax_{\btheta } \e_{\bX}  \{ {\cal L}_n^S(\btheta) \} = \mo$  by the first order condition and concavity, and the $S$-restricted least absolute deviation estimator can be written as
\begin{align} \label{S-LAD}
		\hat{\btheta}_S = \argmax_{\btheta\in \bbr^s} {\cal L}_n^S(\btheta).
\end{align}

We first establish in Lemma~\ref{LAD.concentration} an upper bound for the maximum $\ell_2$-risks of $\hat{\btheta}_S$.

\begin{lemma}  \label{LAD.concentration}
Assume that \eqref{eq3.8} holds and that $\e  |\varepsilon|^\kappa <\infty$ for some $1< \kappa \leq 2$. Then, on the event $\mathcal{E}_0(\tau)$ for $\tau>0$, the sequence of LAD estimators $\{\hat{\btheta}_S: S\subseteq [p], |S| = s\}$ satisfies
\begin{align}  \label{uniform.LAD.bound}
 \max_{S \subseteq [p]: |S| =s }\| \hat{\bSigma}_{SS}^{1/2} \hat{\btheta}_S \|_2  \leq  C_1 \,  a_2^{-1}  \{  s\log(pn) \}^{1/2}   n^{-1/2}
\end{align}
with conditional probability (over the randomness of $\{\varepsilon_i\}_{i=1}^n$) greater than $1- c_1 n^{-1} - c_2  n^{1-\kappa}$, where $C_1, c_1 >0$ are absolute constants and $c_2>0$ is a constant depending only on $a_2$, $\kappa$, $\e |\varepsilon|  $ and $\e  |\varepsilon|^\kappa $.
\end{lemma}

Based on Lemma~\ref{LAD.concentration}, we further study the concentration property of the Wilks expansion for the excess ${\cal L}_n^S(\hat \btheta_S) - {\cal L}_n^S(\mo )$. Since the function ${\cal L}_n^S(\cdot)$ is concave, we use $\nabla  {\cal L}_n^S(\cdot)$ to denote its subgradient. For $\btheta \in \bbr^s$, let $\bzeta^S(\btheta) = {\cal L}_n^S(\btheta) - \e_{\bX}  {\cal L}_n^S(\btheta) $ be the stochastic component of ${\cal L}_n^S(\btheta)$. Then, it is easy to see that
\begin{align}
	\nabla  \bzeta^S(\btheta)  = - 2 \sn  w^S_i(\btheta) \bX_{iS} , \quad \nabla   \e_{\bX}  {\cal L}_n^S(\btheta) = -  \sn \{  2 \P_{\bX}(Y_i \leq \bX_{iS}^\T \btheta) - 1 \}  \bX_{iS},  \label{S-Gradient}
\end{align}
where $w^S_i(\btheta):=  I(Y_i \leq \bX_{iS}^\T \btheta  ) - \P_{\bX}(Y_i  \leq \bX_{iS}^\T \btheta ) $. In particular, we have $\nabla \bzeta^S(\mo ) = - \sn \{ 2 I(\varepsilon_i \leq 0) -  1\} \bX_{iS}$. Recall that $f_\varepsilon$ and $F_\varepsilon$ denote, respectively, the density function and the cumulative distribution function of $\varepsilon$. By the second expression in \eqref{S-Gradient}, $ \nabla \e_{\bX}  {\cal L}_n^S( \btheta ) = - \sn  \{ 2 F_\varepsilon( \bX_{iS}^\T \btheta)-1\} \bX_{iS}$ and
\begin{align}
   \bH_S(\btheta) := - \nabla^2 \e_{\bX}  {\cal L}_n^S( \btheta )  =   2 \sn f_\varepsilon( \bX_{iS}^\T \btheta) \bX_{iS} \bX_{i S}^\T. \label{S-Hessian}
\end{align}
In line with \eqref{infor.mat.sub-model}, we have $\bH^*_S = \bH_S(\mo ) =  2 n f_\varepsilon(0)  \,  \hat{\bSigma}_{SS} $, which is the negative Hessian of $\e_{\bX} {\cal L}_n^S(\mo )$. As in \eqref{normalized.score.sub-model}, define the normalized score
\begin{align}
	\hat{\bxi}_S   =\bH^{* -1/2}_S  \nabla {\cal L}_n^S(\mo) =   \{2n f_\varepsilon(0)\}^{-1/2} \, \hat\bSigma_{SS}^{-1/2} \sn \{ 2 I( \varepsilon_i \leq 0 ) - 1\} \bX_{iS}. \label{S-score}
\end{align}

The following result is a non-asymptotic, conditional version of the Wilks theorem, saying that with high probability, the square root of the excess $\max_{\btheta  } {\cal L}_n^S(\btheta) - {\cal L}_n^S(\mo )$ and the $\ell_2$-norm of the normalized score $\hat{\bxi}_S$ are sufficiently close uniformly over all $s$-subsets $S\subseteq [p]$.

\begin{lemma} \label{Wilks.approximation.2}
Assume that Conditions~\ref{cond3.1} and \ref{cond3.4} are satisfied. Then
\begin{align}
&  \max_{S\subseteq[p]: |S|=s} \Big|  [ 2 \{ {\cal L}_n^S(\hat{\btheta}_S)  - {\cal L}_n^S(\mo) \} ]^{1/2} -  \| \hat{\bxi}_S \|_2 \Big| \nn  \\
 & \qquad \qquad  \leq  C_2  \{ f_{\varepsilon}(0)\}^{-1/2} \big[ \lambda_{\min}^{-1/2}(s)  \{ s\log(pn) \}^{3/2}  n^{-1/2}   + \lambda_{\min}^{-1/4}(s) s\log (pn) n^{-1/4}  \big] \label{wilks.approx.2}
\end{align}
holds with probability greater than $1- c_2   n^{1-\kappa} - c_3  n^{-1}$ whenever $n \geq  C_3 \, \lambda_{\min}^{-1}(s) \{s\log(pn)\}^2$, where $C_2 >0$ is a constant depending only on $a_2, A_2$ and $A_3$, $c_2$ is as in Lemma~\ref{LAD.concentration}, $c_3>0$ is an absolute constant and $C_3 >0$ is a constant depending only on $a_2$ and $A_2$.
\end{lemma}

Further, write $\widetilde{\varepsilon}_i =2 I(\varepsilon_i \leq 0) -1$ and $\widetilde \bX_i = \widetilde{\varepsilon}_i \bX_i$. Note that $\widetilde{\varepsilon}_1, \ldots, \widetilde{\varepsilon}_n$ are i.i.d. Rademacher random variables and thus $\widetilde \bX_1, \ldots, \widetilde \bX_n$ are sub-exponential random vectors. In this notation, we have $\hat{\bxi}_S = \{ 2 n f_\varepsilon(0)  \}^{-1/2}   \sn  \hat{\bSigma}_{SS}^{-1/2}  \widetilde \bX_{iS}$. For each $S\subseteq [p]$, define
$$
\bxi_S  =  \{ 2 n f_\varepsilon(0)  \}^{-1/2}    \sn   {\bSigma}_{SS}^{-1/2} \bX_{iS}.
$$
Then, applying Lemma~\ref{GAR.lem1} with slight modification and the union bound we obtain that, with probability at least $1- c_4  n^{-1}$,
\begin{align}  \label{unif.quadratic.approxi.2}
\max_{S\subseteq [p]: |S| =s } \big| \| \hat{\bxi}_S \|_2 - \| \bxi_S \|_2 \big| \leq  C_4 \{   f_\varepsilon(0)  \}^{-1/2}    s\log(pn)\, n^{-1/2}
\end{align}
for all $n \geq  C_5\, s\log(pn)$, where $c_4>0$ is an absolute constant and $C_4, C_5 >0$ are constants depending only on $A_0$.

Observe that $\e (\wt \bX_i ) = \e  [ \bX_i \{ 2\P(\varepsilon_i\leq 0 |\bX_i )  -1 \} ] = 0$ and $\e (\wt \bX_i \wt \bX_i^\T  ) =\e ( \bX_i  \bX_i^\T  ) = \bSigma$. Hence, it follows from Lemma~\ref{GAR.lem2} that there exists a random variable $T_0 \stackrel{d}{=} R_0(s,p) $ such that for any $\delta\in (0, 1]$,
\begin{align}
	\bigg|   \sqrt{ 2  f_\varepsilon(0)  }   \max_{S\subseteq [p]: |S| =s } \| \bxi_S \|_2 - T_0 \bigg| \leq C_6  \big[ \delta +  \{  s\log(\gamma_s pn ) \}^{1/2} n^{-1/2} +  \{s\log(\gamma_s pn )\}^2 n^{-3/2}  \big] \label{GAR.BE.2}
\end{align}
holds with probability at least $1-C_7 \big[   \delta^{-3}n^{-1/2} \{s \log(\gamma_s pn)\}^2 \, \vee \, \delta^{-4}n^{-1} \{s \log(\gamma_s pn)\}^5 \big] $, where $C_6, C_7 >0$ are constants depending only on $A_0$.

Finally, combining \eqref{wilks.approx.2}, \eqref{unif.quadratic.approxi.2}, \eqref{GAR.BE.2} and \eqref{R*.anti-concentration} proves \eqref{eq3.10}. \qed


\vskip 0.2in


\newpage

\appendix
\section{Appendix A.}
\label{sec7}



In this appendix we prove the technical lemmas appeared in Section~\ref{sec6}.

\subsection{Proof of Lemma~\ref{quadratic.concentration.lemma}}

Define the loss function $\ell(y,z) = y z - b(z)$ for $y,z \in \bbr$. For each $s$-subset $S\subseteq [p]$ and ${\btheta} \in \bbr^s$, define $\bzeta^S({\btheta}) = {\cal L}_n^S({\btheta}) - \e_{\bX} {\cal L}_n^S({\btheta}) = \sn \bzeta^S_i(  \btheta)$, where $\bzeta^S_i(  \btheta) = \ell(Y_i, {\bX}_{i S}^\T {\btheta}) - \e_{\bX}  \ell(Y_i, {\bX}_{i S}^\T {\btheta})$. Note that $\nabla \bzeta^S_i(  \btheta)|_{\btheta = \mo} = \varepsilon_i \bX_{i S}$ with $\varepsilon_i = Y_i-b'(0) $. Thus, we have ${\bV}_0^2 := \var_{\bX}\{ \nabla \bzeta^S( \mo  ) \} = n \phi  b''(0)  \, \hat \bSigma_{SS}$.

For every $\bu \in \bbr^s \setminus \{\mo \}$ and $u \in \bbr$,
\begin{align}
    \e_{\bX} \exp\bigg\{  u \frac{\bu^\T \nabla \bzeta^S( \mo ) }{ \|  {\bV}_0 \bu \|_2} \bigg\}
  & =  \prod_{i=1}^n \e_{\bX} \exp\biggl( u  \frac{  \bu^\T \bX_{i S}  }{ \|  {\bV}_0 \bu \|_2}   \varepsilon_i \biggr) \nn \\
   & =  \prod_{i=1}^n \e_{\bX} \exp\biggl\{  \frac{u}{\sqrt{n} } \times  \frac{  \bu^\T {\bX}_{i S}    }{  ( \bu^\T  \hat \bSigma_{SS}  \bu )^{1/2}} \times \frac{\varepsilon_i}{  ( \var \, \varepsilon_i )^{1/2} } \biggr\}  \nn \\
  &  \leq \exp \biggl \{\frac{1}{2} a_0 u^2   \times  \frac{1}{n} \sn \frac{ (\bu^\T \bX_{i S}  )^2 }{\bu^\T \hat \bSigma_{SS} \bu} \biggr \} = \exp(  a_0  u^2 /2 ). \nn
\end{align}
This verifies condition ($ED_0$) with $\nu_0^2=a_0$ in Theorem~B.3 from the supplement of \cite{SZ2014}. Consequently, taking $ \mathbb{B}^2 = \bH_{S}^{* -1/2} \bV_0^2  \bH_{S}^{* -1/2} = \phi \bI_{s}$ and ${\rm g}=\{  C {\rm tr}(\mathbb{B}^2) \}^{1/2}$ for some $C \geq 2$ there, we have $ \lambda_{\max}(\mathbb{B}^2) = \phi$, ${\rm tr}(\mathbb{B}^2)=\phi s$, ${\rm tr}(\mathbb{B}^4) = \phi^2 s $ and ${\rm x}_c=  \tfrac{1}{2}  ( \tfrac{3 }{2}  C - 1 - \log 3 ) s  \geq \frac{3}{4}(C-2) s$. This implies that almost surely on the event $\{\hat{\bSigma}_{SS} \succ \mo \}$, with conditional probability at least  $1-2e^{-t} - 8.4 \, e^{-{\rm x}_c}$,
\begin{align*}
 \| \hat \bxi_S \|_2^2 \leq a_0 \phi \times \begin{cases}
 s +  ( 8ts )^{1/2} , \quad  &{\rm if}~ 0\leq t\leq \frac{1}{18} (2s)^{1/2}, \\
 	 s+ 6t , \quad &  {\rm if}~ \frac{1}{18} (2s)^{1/2} < t \leq {\rm x}_c.
 \end{cases}
\end{align*}
Finally, letting $C\to \infty$ proves \eqref{quadratic.concentration.ineq}. \qed

\subsection{Proof of Lemma~\ref{Wilks.approximation}}

We prove this lemma by applying the conditional version of Theorem~2.3 in \cite{S2013}. To this end, we need to verify conditions ($ED_0$), ($ED_2$), ($\mathcal{L}_0$), ($\mathcal{I}$) and ($\mathcal{L}$). In line with the notation used therein, we fix $S\subseteq [p]$ and write
$$
	\bD^2(\btheta) = -\nabla^2 \e_{\bX} \{ {\cal L}_n^S(\btheta) \} = \sn b''(\bX_{iS}^\T \btheta) \bX_{iS} \bX_{iS}^\T, \quad \bD_0^2 = \bD^2(\mo) = n b''(0) \, \hat{\bSigma}_{SS} .
$$

The validity of  ($ED_0$) is guaranteed from the proof of Lemma~\ref{quadratic.concentration.lemma}, and ($ED_2$) is automatically satisfied with $\omega \equiv 0$ since $\nabla^2 \bzeta^S(\btheta)$ vanishes for all $\btheta \in \bbr^s$. Turning to ($\mathcal{L}_0$), observe that
\begin{align}
	&	\big\|  \bD_0^{-1}  \bD^2(\btheta) \bD_0^{-1} - \bI_s \big\|  \nn \\
	& =  \bigg\|\bD_0^{-1}  \sn \{b''(\bX_{i S}^\T \btheta) - b''(0)\} \bX_{i S} \bX_{i S}^\T \bD_0^{-1} \bigg\| \nn \\
	& = \bigg\| \bD_0^{-1}  \sn  b'''( \eta_i)  \bX_{i S}^\T \btheta \,  \bX_{i S} \bX_{i S}^\T \bD_0^{-1}  \bigg\| ,  \label{L0.1}
\end{align}
where $\eta_i$ lies between $0$ and $\bX_{i S}^\T \btheta$. For $r>0$, define $\Theta_0(r) = \{  \btheta  \in \bbr^s :  \|   \bD_0  \btheta    \|_2 \leq r  \}$. On the event $\mathcal{E}_0(\tau)$ for some $\tau>0$ and for $\btheta \in \Theta_0(r)$,
\begin{align}
& | \bX_{i S}^\T \btheta |   = |  \btheta^\T  \bD_0 \bD_0^{-1} \bX_{i S}|   \leq   \| \bD_0^{-1} \bX_{i S} \|_2 \leq    \{ n b''(0)\}^{-1/2}  \tau^{1/2} r .  \label{perturbation}
\end{align}
This together with \eqref{L0.1} implies that
\begin{align}
	\big\| \bD_0^{-1} \bD^2(\btheta) \bD_0^{-1} - \bI_s \big\|   \leq   \frac{   \max_{|t | \leq    \{n b''(0)\}^{-1/2}   \tau^{1/2} r }|b'''(t)| }{  \{ b''(0) \}^{3/2} } \frac{\tau^{1/2} r }{n^{1/2}}  := \delta(\tau, r)  .  \label{L0.2}
\end{align}

Recalling that $\bV_0^2 = \var_{\bX}\{ \bzeta^S(\mo) \} = \phi \bD_0^2$, ($\mathcal{I}$) is satisfied with $a =\phi^{1/2}$.

To verify ($\mathcal{L}{{\rm r}}$), define $g(t)= b'(0)t -b(t)$ so that $g'(t)= b'(0)-b'(t)$ and $g''(t)=-b''(t)$. Then, for any $\btheta \in \bbr^s$ satisfying $\| \bD_0 \btheta \|_2 = r >0$, it follows from the second-order Taylor expansion that
\begin{align}
	 & -2\{ \e_{\bX} {\cal L}_n^S( \btheta )  - \e_{\bX} {\cal L}_n^S( \mo  )  \} = -2 \sn \{  g(\bX_{iS}^\T \btheta) - g(0) \} \nn \\
	& =  -2 \sn \big\{ g'(0) \bX_{iS}^\T \btheta + \tfrac{1}{2} g''(\eta_i)(\bX_{iS}^\T \btheta )^2 \big\} =  \sn   b''( \eta_i)(\bX_{i S}^\T \btheta )^2 , \label{Lr.1}
\end{align}
where $\eta_i$ is a point lying between $0$ and $\bX_{i S}^\T \btheta$. On the event $\mathcal{E}_0(\tau)$, the right-hand side of \eqref{Lr.1} is further bounded from below by
\begin{align}
  r^2 \{ b''(0)\}^{-1}\min_{|t| \leq \{n b''(0)\}^{-1/2} \tau^{1/2} r } b''(t)  . \nn
\end{align}
When $  \| \bD_0 \btheta \|_2 = r  \leq  \{ n b''(0)/\tau \}^{1/2}$, $ -2\{ \e_{\bX} {\cal L}_n^S( \btheta )  - \e_{\bX} {\cal L}_n^S( \mo  )  \} $ is bounded from below by $a_1 r^2 $ for $a_1$ as in \eqref{eq3.5}. Further, from the convexity of the function $\btheta \mapsto - \e_{\bX} \{{\cal L}_n^S(\btheta) - {\cal L}_n^S(\mo) \} $, we see that  $- \e_{\bX} \{{\cal L}_n^S(\btheta) - {\cal L}_n^S(\mo) \} \geq    a_1 r \{ n b''(0)/\tau \}^{1/2}$, for all $\btheta$ satisfying $\| \bD_0 \btheta \|_2 =r \geq  \{ n b''(0)/\tau \}^{1/2}$. Define the function $r\mapsto b(r)$ as
\begin{align}
	b(r) = \begin{cases}
		a_1 \quad & \mbox{ if }  0\leq r\leq  \{ n b''(0)/\tau \}^{1/2} , \\
	    a_1 r^{-1} \{ n b''(0)/\tau \}^{1/2}  \quad & \mbox{ if }  r> \{ n b''(0)/\tau \}^{1/2} .
	\end{cases}  \label{def.br}
\end{align}
By definition, $rb(r)$ is non-decreasing in $r\geq 0$ and for $\btheta \in \bbr^s$ satisfying $\| \bD_0 \btheta \|_2 =r$,
\begin{align}
 	\frac{ -2 \e_{\bX} \{{\cal L}_n^S(\btheta) - {\cal L}_n^S(\mo) \}  }{\| \bD_0 \btheta \|_2^2 } \geq b(r).   \label{global.id}
\end{align}

With the above preparations, we apply Theorem~2.3 in \cite{S2013} with slight modification on the constant. In view of \eqref{Delta.st} and \eqref{def.br}, set
\begin{align}
	r_0 = 2(\phi a_0)^{1/2} a_1^{-1} \big[ s + 6\big(  s\log \tfrac{ep}{s} + \log n \big)  \big]^{1/2} ,  \label{def.r0.GLIM}
\end{align}
such that Condition~2.3 there is satisfied on $\mathcal{E}_0(\tau)$ whenever $n\geq \{b''(0)\}^{-1} r_0^2 \tau$. Hence, it follows from Theorem~2.3 in \cite{S2013} and the union bound that, conditional on the event $\mathcal{E}_0(\tau)$,
\begin{align}
	\P_{\bX} \bigg(  \max_{S\subseteq [p]: |S|=s} \big|  [  2\{ {\cal L}_n^S(\hat{\btheta}_S) - {\cal L}_n^S(\mo) \}  ]^{1/2}  - \| \hat{\bxi}_S \|_2 \big| \leq 5   \delta(\tau, r_0) r_0 \bigg) \leq 5 n^{-1} ,
\end{align}
where $\delta(\tau, r)$ and $r_0$ are as in \eqref{L0.2} and \eqref{def.r0.GLIM}, respectively. This proves \eqref{Wilks.concentration} by properly choosing $C_1$ and $C_2$.	\qed

\subsection{Proof of Lemma~\ref{GAR.lem1}}

To begin with, note that for each $s$-subset $S\subseteq [p]$, $\bZ_{1S}, \ldots, \bZ_{nS}$ are i.i.d. $s$-dimensional random vectors with mean zero and covariance matrix $\phi  \bSigma_{SS}$. By \eqref{normalized.score.sub-model} and \eqref{xiS.def},
\begin{align}
  \| \hat{\bxi}_S \|_2^2 -   \|  {\bxi}_S \|_2^2   = \bxi_S^\T  \big( \bSigma_{SS}^{ 1/2} \hat \bSigma_{SS}^{-1} \bSigma_{SS}^{1/2}  - \bI_s \big)  \bxi_S. \nn
\end{align}
Write $\bbX_S = (\bX_{1S}, \ldots, \bX_{nS})^\T \in \bbr^{n\times s}$, then $\bbX_S \bSigma_{SS}^{-1/2}$ is an $n\times s$ matrix whose rows are independent sub-Gaussian random vectors in $\bbr^s$. Further, observe that $\bX_{iS}= \bP_S \bX_i$ and $\bSigma_{SS} = {\bP}_S \bSigma {\bP}_S^\T$, where $ \bP_S \in \bbr^{s\times p}$ is a projection matrix. Under Condition~\ref{cond3.1}, $\| \bu^\T \bSigma_{SS}^{-1/2} \bX_{iS} \|_{\psi_2} = \| \bu^\T \bSigma_{SS}^{-1/2}  \bP_S  \bSigma_{SS}^{1/2}  \bU   \|_{\psi_2} \leq A_0 \|  \bSigma_{SS}^{1/2}  \bP^\T_S \bSigma_{SS}^{-1/2}\bu \|_2 = A_0$ for $\bu\in\mathbb{S}^{s-1}$. Then, it follows from \eqref{hatS.positivity} that for all sufficient large $n$ so that $\delta\leq \frac{1}{2}$, $\|  \bSigma_{SS}^{ 1/2} \hat \bSigma_{SS}^{-1} \bSigma_{SS}^{ 1/2} - \bI_s \| \leq 2\delta$ and hence,
\begin{align}
  &  |  \|  \hat{\bxi}_S \|_2  -   \|   {\bxi}_S \|_2   | = \frac{ |  \|  \hat{\bxi}_S \|_2^2  -   \|   {\bxi}_S \|^2_2   | }{  \|  \hat{\bxi}_S \|_2 +  \|   {\bxi}_S \|_2 } \nn \\
  & \leq \| \bxi_S \|_2^{-1} \times  |   \|  \hat{\bxi}_S \|_2^2 -   \|   {\bxi}_S \|_2^2  | \leq  2C_3 (s \vee t)^{1/2} n^{-1/2}     \times   \|  \bxi_S \|_2  .	\label{quadratic.diff}
\end{align}

Next we upper bound the quadratic term $  \| \bxi_S \|_2$. First we show that $\bSigma_{SS}^{-1/2} \bZ_{iS} = \phi^{1/2} \, \bSigma_{SS}^{-1/2}\wt{\varepsilon}_i \bX_i$ are sub-exponential random vectors, where $\wt{\varepsilon}_i  := \varepsilon_i/ ( \var\,\varepsilon_i)^{1/2}$. In fact, for every $\bu \in \mathbb{S}^{s-1}$, $\|  \bu^\T \bSigma_{SS}^{-1/2} \bZ_{iS} \|_{\psi_1} \leq  2 \| \wt \varepsilon_i \|_{\psi_2} \| \bu^\T \bSigma_{SS}^{-1/2} \bX_{iS} \|_{\psi_2} \leq 2 A'_0 A_0$, where $A'_0>0$ is a constant depending only on $a_0$ in Condition~\ref{cond3.1}. Following the proof of Lemma~5.15 in \cite{V2012}, we derive that for every $\bu \in \bbr^s$ satisfying $\| \bu \|_2\leq  \phi^{-1/2} ( 4e A'_0 A_0 )^{-1 } \sqrt{n} $,
\begin{align}
	\log \e \exp( \bu^\T \bxi_S) & = \sn \log \e \exp( n^{-1/2} \bu^\T \bSigma_{SS}^{-1/2} \bZ_{iS} ) \nn \\
	& \leq 2 e^2 \|\bu\|_2^2 \,  n^{-1} \sn \big\| (\bu/\|\bu\|_2)^\T \bSigma_{SS}^{-1/2} \bZ_{iS}  \big\|_{\psi_1}^2  \nn \\
	& \leq  ( 4e  A'_0 A_0 )^2 \phi  \frac{\| \bu \|_2^2 }{2}. \nn
\end{align}
Consequently, applying Corollary~1.12 in the supplement of \cite{S2012} with ${\rm g}= \sqrt{n}$, $\mathbb{B}=\bI_s$ and ${\rm x}_c=\frac{3}{4} n - \frac{1}{2}(1+\log 3)s \geq \frac{3}{4} n -\frac{3}{2}s$ to the random vector $(4e A'_0 A_0)^{-1}  \phi^{-1/2} \bxi_S$ yields that, for every $0\leq t\leq {\rm x}_c$,
\begin{align}
	\P\big[ \|  \bxi_S \|_2  \geq  4 e A'_0 A_0 \{  \phi \Delta(s,t) \}^{1/2} \big] \leq 2 e^{-t} + 8.4 \,e^{-{\rm x}_c}. \label{quadratic.deviation}
\end{align}

Finally, combining \eqref{quadratic.diff} and \eqref{quadratic.deviation} completes the proof of \eqref{quadratic.approxi}. \qed

\subsection{Proof of Lemma~\ref{GAR.lem2}}

First, observe that
\begin{align}
	\max_{S\subseteq[p]: |S|=s} \| \bxi_S \|_2
= \max_{\bu \in \mathcal{F}(s,p)}  n^{-1/2} \sn \frac{\bu^\T \bZ_i }{ ( \bu^\T \bSigma \bu )^{1/2}}, \nn
\end{align}
where $\mathcal{F}(s,p)  = \{  \bx \mapsto  \bu^\T \bx   : \bu  \in \mathbb{S}^{p-1} ,  \|  \bu \|_0 \leq s  \}$. Recall that $\bZ_1,\ldots, \bZ_n$ are i.i.d. $p$-dimensional centered random vectors with covariance matrix $\e (\bZ_i \bZ_i^\T )= \phi \bSigma$. As in the proof of Lemma~\ref{GAR.lem1}, we have for any $\bu\in \mathbb{S}^{p-1}$,
\begin{align*}
	\| \phi^{-1/2} \bu^{{\rm T}} \bZ_{i} \|_{\psi_1}  \leq 2 \| \varepsilon_i / (\var\, \varepsilon_i)^{1/2} \|_{\psi_2} \, \| \bu^\T \bSigma^{1/2} \bU_i \|_{\psi_2} \leq 2 A'_0 A_0 ( \bu^\T \bSigma \bu )^{1/2} .
\end{align*}
Consequently, it follows from Lemma~7.5 in \cite{FSZ2015} that there exists a random variable $T_0 \stackrel{d}{=} R_0(s,p) = \max_{\bu \in \mathcal{F}(s,p)} \frac{ \bu^{{\rm T}}\bG }{ ( \bu^{{\rm T}} \bSigma \bu )^{1/2} }$ for $\bG \sim N(\mo, \bSigma)$  such that, for any $\delta \in (0, 1]$,
\begin{align}
	& \P\bigg\{ \bigg|  \max_{S\subseteq[p]: |S|=s} \| \phi^{-1/2} \bxi_S \|_2 -  T_0 \bigg|  \geq C_1  A'_0 A_0 \bigg(  \delta +  \frac{\gamma^{1/2}_{s,p,n}}{\sqrt{n}}    +  \frac{\gamma_{s,p,n}^2 }{n^{3/2}} \bigg) \bigg\} \nn \\
	& \qquad \qquad \qquad \qquad  \qquad \qquad \qquad  \qquad \leq C_2 \bigg[    \frac{\{s \log(\gamma_s pn)\}^2}{\delta^3 \sqrt{n}} +   \frac{\{s\log(\gamma_s pn )\}^5}{\delta^4 n} \bigg] , \nn
\end{align}
where $\gamma_{s,p,n} = s\log \frac{\gamma_s ep}{s} + \log n$ and $C_1, C_2 >0$ are absolute constants. This proves \eqref{GAR.BE}. \qed

\subsection{Proof of Lemma~\ref{LAD.concentration}} The proof employs techniques from empirical process theory which modify the arguments used in \cite{W2013}. To begin with, note that
$$
	\hat{\btheta}_S =  \arg\min_{\btheta \in \bbr^s }  f(\btheta) : = \arg\min_{\btheta \in \bbr^s } \| \bY - \bbX_S \btheta \|_1.
$$
Under the null model, $\bY= \bbX \bbeta^* +\beps = \bbX_S \btheta^* + \beps$ with $\btheta^* = \mo$. Then the sub-differential of $f(\btheta)$ at $\btheta =\mo $ can be written as $\nabla f(\mo) = - \bbX_S^\T \sgn(\beps)$, where $\sgn(\beps)=(\sgn (\varepsilon_1), \ldots,  \sgn (\varepsilon_n))^\T$ with $\sgn(u) :=I(u >0)-I(u<0)$. Define $\bz  = (z_1, \ldots , z_n)^\T = \sgn(\beps)$, and note that $z_1,\ldots, z_n$ are i.i.d. random variables satisfying $\P(z_i=1)=\P(z_i=-1)=1/2$.

Since $\hat{\btheta}_S$ minimizes $\| \bY - \bbX_S \btheta \|_1$ over $\bbr^s$, we have the following basic inequality
\begin{align}
	\| \bY -  \bbX_S \hat{\btheta}_S  \|_1 = \| \bbX_S \hat{\btheta}_S - \beps \|_1    \leq \| \beps \|_1. \label{basic.inequality}
\end{align}
Further, define a random process $\{Q(\btheta)\}$ indexed by $\btheta \in \bbr^s$:
\begin{align}
	Q(\btheta) =  n^{-1/2} \sn \big(  | \bX_{iS}^\T \btheta - \varepsilon_i | - |\varepsilon_i|  \big). \label{Q.theta}
\end{align}
In what follows, we prove that with overwhelmingly high probability , $Q( {\btheta})$ is concentrated around its expectation $Q_{\bbX}(\btheta) := \e_{\bX} \{ Q(\btheta)\}$ uniformly over $\btheta \in \bbr^s$ via a straightforward adaptation of the peeling argument.

For $\delta_1 >0$ and $\ell=1, 2, \ldots$, consider the following sequence of events
\begin{align}
\mathcal{G}(\delta_1 ) =\big\{ \btheta \in \bbr^s :  \| \hat{\bSigma}_{SS}^{1/2} \btheta \|_2   \geq \delta_1 \big\}, \quad \mathcal{G}_\ell(\delta_1) =\big\{ \btheta \in \bbr^s : \alpha^{\ell-1}\delta_1 \leq \| \hat{\bSigma}_{SS}^{1/2} \btheta \|_2   \leq \alpha^\ell \delta_1 \big\}, \label{event.G}
\end{align}
where $\alpha= \sqrt{2}$. Here, $\delta_1$ can be regarded as a tolerance parameter, and it is easy to see that $\mathcal{G}(\delta_1)= \cup_{\ell=1}^\infty \mathcal{G}_\ell(\delta_1)$. For $R>0$, set $\mathcal{V}(R) =  \{ \btheta \in \mathcal{G}(\delta_1) :   \| \hat{\bSigma}_{SS}^{1/2}  \btheta \|_2  \leq R\}$ and let $\Delta(R)$ be the maximum deviation over the elliptic vicinity $\mathcal{V}(R)$:
\begin{align}
\Delta(R) =  \max_{\btheta \in \mathcal{V}(R) } | Q(\btheta) - Q_{\bbX}(\btheta) |.
\end{align}
For every $\btheta \in \bbr^s$, define the rescaled vector $\wt \btheta = \hat \bSigma_{SS}^{1/2} \btheta$ such that
$$
	\Delta(R) =  \max_{ \delta_1 \leq \| \wt \btheta \|_2 \leq R } \Big| Q(\hat \bSigma_{SS}^{-1/2} \wt \btheta) - Q_{\bbX}(\hat \bSigma_{SS}^{-1/2} \wt \btheta) \Big| .
$$
For every $0<\epsilon \leq R$, there exists an $\epsilon$-net $\mathcal{N}_\epsilon$ of the Euclidean ball $\mathbb{B}_2^s(R)$ with cardinality bounded by $ \big(1+\frac{2R}{\epsilon} \big)^s$. For $\wt \btheta_1, \wt \btheta_2 \in \mathbb{B}^s_2(R)$ satisfying $\| \wt \btheta_1 - \wt \btheta_2 \|_2 \leq  \epsilon$, observe that
\begin{align}
	\Big| Q( \hat{\bSigma}_{SS}^{-1/2} \wt \btheta_1) - Q(\hat{\bSigma}_{SS}^{-1/2} \wt \btheta_2) \Big| &   \leq   n^{-1/2} \sn \Big| \bX_{iS}^\T \hat{\bSigma}_{SS}^{-1/2} ( \wt \btheta_1 - \wt \btheta_2) \Big| \nn \\
	& \leq  \Big\| \bbX_S \hat{\bSigma}_{SS}^{-1/2} (\wt \btheta_1 - \wt \btheta_2) \Big\|_2 \leq  \epsilon   n^{1/2}   . \nn
\end{align}
Then, it is easy to see that
\begin{align}
	 \Delta(R)  \leq \max_{\wt \btheta \in \mathcal{N}_\epsilon } \Big| Q(\hat{\bSigma}_{SS}^{-1/2} \wt \btheta) - Q_{\bbX}(\hat{\bSigma}_{SS}^{-1/2} \wt \btheta) \Big|   + 2\epsilon  n^{1/2}  . \label{Delta.discretization}
\end{align}
For each $\wt \btheta \in  \mathbb{B}^s_2(R)$ fixed, $Q(\hat{\bSigma}_{SS}^{-1/2} \wt \btheta) - Q_{\bbX}(\hat{\bSigma}_{SS}^{-1/2} \wt \btheta)$ is a sum of independent random variables with zero means and for $i=1,\ldots, n$, $ | | \bX_{iS}^\T  \hat{\bSigma}_{SS}^{-1/2} \wt \btheta - \varepsilon_i | - |\varepsilon_i|  | \leq | \bX_{iS}^\T \hat{\bSigma}_{SS}^{-1/2} \wt  \btheta|$. Therefore, it follows from Hoeffding's inequality that for every $t>0$,
\begin{align}
	& \P_{\bX} \Big\{ \Big| Q(\hat{\bSigma}_{SS}^{-1/2} \wt \btheta) - Q_{\bbX}(\hat{\bSigma}_{SS}^{-1/2} \wt \btheta) \Big| \geq t \Big\} \nn \\
	& \leq 2\exp\left\{ -\frac{ n t^2}{2 \sn (\bX_{iS}^{{\rm T}} \hat{\bSigma}_{SS}^{-1/2} \wt \btheta )^2  } \right\} = 2\exp\bigg( -\frac{   t^2}{2  \| \wt \btheta \|_2^2   }  \bigg)  .\nn
\end{align}
In other words, for every $\wt \btheta \in  \mathbb{B}^s_2(R)$ and $\delta  > 0$,
\begin{align}
	\Big| Q(\hat{\bSigma}_{SS}^{-1/2} \wt \btheta) - Q_{\bbX}(\hat{\bSigma}_{SS}^{-1/2} \wt \btheta) \Big|  \leq  ( 2\delta )^{1/2}  \| \wt  \btheta  \|_2 \leq  ( 2\delta )^{1/2}  R   \nn
\end{align}
holds with probability at least $1-2e^{-\delta}$. This, together with the union bound yields
\begin{align}
 \P_{\bX} \bigg\{ \max_{\wt \btheta \in \mathcal{N}_\epsilon } \Big| Q(\hat{\bSigma}_{SS}^{-1/2} \wt \btheta) - Q_{\bbX}(\hat{\bSigma}_{SS}^{-1/2} \wt \btheta) \Big|  \geq  ( 2\delta )^{1/2} R  \bigg\}
 \leq \exp\big\{ s\log\big( 1+\tfrac{2R}{\epsilon} \big) - \delta \big\} .   \label{discrete.Delta.tail.prob}
\end{align}
In particular, by taking $\epsilon=R n^{-1}$ in \eqref{Delta.discretization} and $\delta=s\log ( 1+\tfrac{2R}{\epsilon} ) + t \leq 2s\log n+ t$ in \eqref{discrete.Delta.tail.prob} we conclude that
\begin{align}
	\P_{\bX} \big\{ \Delta(R) \geq  R (2t)^{1/2} + 2R ( s\log n )^{1/2} + 2 R n^{-1/2} \big\} \leq 2e^{-t} \label{cont.Delta.tail.prob}
\end{align}
holds almost surely on the event $\mathcal{E}_0(\tau)$ for any $\tau>0$.

In particular, by taking $t=c n R^2 $ in \eqref{cont.Delta.tail.prob} for some $c>0$ to be specified below \eqref{concentration.LAD} and the union bound, we have
\begin{align}
	& \P_{\bX} \Big[  \exists \, \btheta \in  \mathcal{G}(\delta_1), {\rm s.t.} \,  | Q(\btheta) - Q_{\bbX}(\btheta) | \geq  2^{3/2} \| \wt \btheta \|_2 \big\{   \| \wt \btheta \|_2 (cn )^{1/2}    +  (s \log n)^{1/2} + n^{-1/2} \big\} \Big] \nn \\
	& \leq  \sum_{\ell=1}^\infty   \P_{\bX} \Big[  \exists \, \btheta \in  \mathcal{G}_\ell(\delta_1), {\rm s.t.} \,  | Q(\btheta) - Q_{\bbX}(\btheta) | \geq      (\alpha^\ell \delta_1)^2 (2c n )^{1/2} + 2\alpha^\ell \delta_1  \big\{ (s\log n)^{1/2} + n^{-1/2} \big\} \Big] \nn \\
	& \leq   \sum_{\ell=1}^\infty \P_{\bX}\Big[  \Delta(\alpha^\ell \delta_1) \geq (\alpha^\ell \delta_1)^2 (2c n )^{1/2}  + 2\alpha^\ell \delta_1 \big\{ (s\log n)^{1/2} + n^{-1/2} \big\} \Big]  \nn \\
	& \leq 2 \sum_{\ell=1}^\infty \exp\{  - c n(\alpha^\ell \delta_1)^2\}  \leq 2 \sum_{\ell=1}^\infty  \exp\{  -2 c \ell   \log(\alpha) n \delta_1^2 \} \leq \frac{2\exp(-c_0 n \delta_1^2)}{1- \exp(-c_0 n \delta_1^2)}, \nn
\end{align}
where $c_0 = c \log 2$. This implies that with probability at least $1-4\exp(-c_0 n \delta_1^2)$,
\begin{align}  \label{concentration.Q}
	 | Q(\btheta) - Q_{\bbX}(\btheta) | \leq  2^{3/2} \sqrt{c} \, \| \hat{\bSigma}_{SS}^{1/2}\btheta \|_2^2   + 2^{3/2} \| \hat{\bSigma}_{SS}^{1/2}\btheta \|_2 \big\{  (s \log n)^{1/2} + n^{-1/2} \big\}
\end{align}
holds for all $\btheta \in \mathcal{G}(\delta_1)$ whenever $n\geq c^{-1} \delta_1^{-2} $.

For the (conditional) expectation
$$
	Q_{\bbX}(\btheta) = n^{-1/2} \sn \e_{\bX} \big( |\bX_{iS}^\T \btheta - \varepsilon_i | - |\varepsilon_i| \big) = n^{-1/2} \big(  \e_{\bX} \| \bbX_{S}\btheta - \beps \|_1 - \e\|\beps \|_1 \big),
$$
applying Lemmas~5 and 6 in \cite{W2013} with slight modifications gives
\begin{align}  \label{mean.Q.lower.bound}
Q_{\bbX}(\btheta) \geq \begin{cases}
 \frac{1}{4\sqrt{n}} \| \bbX_S \btheta \|_1  =   \frac{ \sqrt{n} }{4 } \| n^{-1}\bbX_S \btheta \|_1 \ \ & \mbox{ if }  \| \bbX_S \btheta \|_1 \geq \frac{2n}{a_2}, \\
 \frac{a_2}{8 \sqrt{n} } \| \bbX_S \btheta \|_2^2 = \frac{a_2 \sqrt{n} }{8} \| \hat{\bSigma}_{SS}^{1/2} \btheta \|_2^2  \ \ & \mbox{ if }  \| \bbX_S \btheta \|_1 < \frac{2n}{a_2},
\end{cases}
\end{align}
where $a_2$ is as in Condition~\ref{cond3.4}. For the sequence of LAD estimators $\{ \hat \btheta_S : S\subseteq [p], |S| =s \}$, from \eqref{basic.inequality} it can be seen that $\| \bbX_S \hat{\btheta}_S \|_1  \leq  \| \bbX_S \hat{\btheta}_S - \beps \|_1 +  \| \beps \|_1 \leq 2\|\beps \|_1$, and hence
$$
	\max_{S\subseteq[p]: |S|=s}  \| n^{-1} \bbX_S \hat{\btheta}_S \|_1  \leq  2  \bigg\{ \e |\varepsilon|+ n^{-1} \sn( |\varepsilon_i|-\e |\varepsilon_i| ) \bigg\}.
$$
For every $t>0$ and $1< \kappa \leq 2$, by Markov's inequality we have
\begin{align}
	\P\bigg\{  \sn( |\varepsilon_i|-\e |\varepsilon_i| )  \geq t \bigg\} \leq t^{-\kappa} \, \e \bigg| \sn( |\varepsilon_i|- \e |\varepsilon_i| )  \bigg|^\kappa \leq 4^{2-\kappa} t^{-\kappa} n \, \e |\varepsilon|^\kappa  ,\nn
\end{align}
where we used the inequality $|1+x|^\kappa \leq 1+ \kappa x + 2^{2-\kappa}|x|^\kappa$ for $1< \kappa \leq 2$ and $x\in \bbr$. The last two displays together imply that, with probability at least $1-\delta_2$,
\begin{align}
 \max_{S\subseteq[p]: |S|=s}  \| n^{-1} \bbX_S \hat{\btheta}_S \|_1 \leq 2 \e  |\varepsilon |   \big\{ 1  + 4^{(2-\kappa)/\kappa}  ( \e |\varepsilon| )^{-1}  (   \e |\varepsilon|^\kappa  )^{1/\kappa}  \delta_2^{-1/\kappa}  n^{-1+ 1/\kappa} \big\} .   \nn
\end{align}
By Condition~\ref{cond3.4}, we have $ a_2 \e  |\varepsilon |  <  1$. Therefore, as long as the sample size $n$ satisfies
\begin{align}
	 n \geq    \bigg\{ \frac{ 4^{2-q }a_2^\kappa \, \e |\varepsilon|^\kappa   }{( 1 - a_2 \, \e |\varepsilon| )^\kappa} \bigg\}^{1/(\kappa - 1)}     \delta_2^{-1/(\kappa -1)} , \label{n.constraint.1}
\end{align}
the event
\begin{align}  \label{event2.def}
  \mathcal{E}_1  := \bigg\{ \max_{S\subseteq[p]: |S|=s}  \| n^{-1} \bbX_S \hat{\btheta}_S \|_1 \leq 2 a_2^{-1} \bigg\}
\end{align}
occurs with probability at least $1-\delta_2$.

Now, by \eqref{basic.inequality}, we have $  Q(\hat{\btheta}_S) \leq 0$ and thus $ -\{  Q(\hat{\btheta}_S) -  Q_{\bbX}(\hat{\btheta}_S) \} \geq   Q_{\bbX}(\hat{\btheta}_S)$ holds for every $s$-subset $S\subseteq [p]$. Together with \eqref{concentration.Q}--\eqref{event2.def} and the union bound, this implies that on the event $\mathcal{E}_0(\tau) \cap \mathcal{E}_1$ for any $\tau >0$,
\begin{align}
	\max_{S\subseteq [p]: |S| = s} \| \hat{\bSigma}_{SS}^{1/2} \hat{\btheta}_S \|_2 \leq \min\bigg[ \delta_1 ,   32 \sqrt{2} \, a_2^{-1}   \bigg\{ \bigg( \frac{s\log n}{n} \bigg)^{1/2} + \frac{1}{n} \bigg\} \bigg]  \label{concentration.LAD}
\end{align}
holds with (conditional) probability $1-4 {p \choose s}\exp( - c_0 n \delta_1^2) - \delta_2$, provided that the sample size $n$ satisfies $n\geq 2 \cdot 32^2  ( a_2 \delta_1 )^{-2}$ and \eqref{n.constraint.1}.

Finally, taking
$$
	\delta_1 = \frac{32}{a_2} \sqrt{\frac{2}{\log(2)}}   \bigg(  \frac{s\log \frac{ep}{s} + \log n}{ n} \bigg)^{1/2}  \ \ \mbox{ and } \ \  \delta_2  =   \frac{4^{2-q} a_2^\kappa \, \e |\varepsilon|^\kappa}{(1-a_2 \, \e |\varepsilon |)^\kappa} \frac{1}{n^{\kappa -1}}
$$
in \eqref{concentration.LAD} proves \eqref{uniform.LAD.bound}. \qed

\subsection{Proof of Lemma~\ref{Wilks.approximation.2}}  We prove this lemma by employing the arguments similar to those used in \cite{S2013}, where the likelihood function ${\cal L}(\btheta)$ is assumed to be twice differentiable with respect to $\btheta$. It is worth noticing that both Conditions~($\mathcal{L}$) and ($ED_2$) in \cite{S2013} are not satisfied in the current situation. We provide here a self-contained proof in which Lemma~\ref{LAD.concentration} also plays an important role.

\medskip
\noindent
{\it \textbf{Step~1: Local linear approximation of $\nabla {\cal L}^S_{n }(\btheta)$}}. Let $\chi_1^S(\btheta)$ be the normalized residual of the local linear approximation of $\nabla {\cal L}^S_{n }(\btheta)$ given by
\begin{align}
	\chi_1^S(\btheta) &=   \bD_0^{-1}\{  \nabla {\cal L}_n^S(\btheta) - \nabla {\cal L}_n^S(\mo) + \bD_0^2\btheta \} \nn \\
	& =   \bD_0^{-1}\{ \bU(\btheta )   +   \nabla \e_{\bX} {\cal L}_n^S(\btheta) -   \nabla \e_{\bX} {\cal L}_n^S(\mo) + \bD_0^2\btheta \} , \label{def.chi1}
\end{align}
where $\bU(\btheta )  = \nabla \bzeta^S(\btheta) - \nabla \bzeta^S(\mo)$ and $\bD_0^2 = - \nabla^2 \e_{\bX} \{ {\cal L}^S_n( \mo) \} = 2f_\varepsilon(0) \sn \bX_{iS} \bX_{iS}^\T$. Then it follows from the mean value theorem that
\begin{align}
	\e_{\bX} \{ \chi^S_1(\btheta)  \} = \{ \bI_s - \bD_0^{-1} \bD^2(\wt \btheta )  \bD_0^{-1} \} \bD_0  \btheta ,  \label{mean.chi-theta}
\end{align}
where $\bD^2(\btheta) = - \nabla^2 \e_{\bX} \{  {\cal L}_n^S(\btheta) \} = 2\sn f_\varepsilon( \bX_{iS}^\T \btheta ) \bX_{iS} \bX_{iS}^\T$ and $\wt \btheta =   \lambda \btheta $ for some $0\leq  \lambda  \leq 1$. As before, for every $r\geq 0$, define the local elliptic neighborhood of $\mo$ as
\begin{align}
	\Theta_0(r) = \{ \btheta \in \bbr^s : \| \bD_0   \btheta \|_2 \leq r \}. \nn
\end{align}
On the event $\mathcal{E}_0(\tau )$ for some $\tau >0$,
\begin{align}
  | \bX_{iS}^\T \btheta | \leq \| \bD_0 \btheta \|_2 \| \bD_0^{-1} \bX_{iS} \|_2  \leq \{ 2n f_\varepsilon(0) \}^{-1/2} \tau^{1/2} r
\end{align}
for all $\btheta \in \Theta_0(r)$. Thus it follows from the Taylor expansion that for $r \leq  \{ 2nf_\varepsilon(0) /\tau \}^{1/2}$,
\begin{align}
  & \big\| \bI_s - \bD_0^{-1} \bD^2(\wt \btheta ) \bD_0^{-1} \big\|  \nn \\
  & = 2 \bigg\| \bD_0^{-1}  \sn \{ f_\varepsilon( \bX_{iS}^\T \wt \btheta ) - f_\varepsilon (0) \} \bX_{iS } \bX_{iS}^\T   \bD_0^{-1}   \bigg\|  \leq \frac{A_3  }{ \sqrt{2} f^{3/2}_\varepsilon(0)}  \frac{ \tau^{1/2} r }{ n^{1/2} }   := \delta( \tau, r ) . \label{def.delta-r}
\end{align}
Together, \eqref{mean.chi-theta} and \eqref{def.delta-r} imply that under the same constraint for \eqref{def.delta-r},
\begin{align}
	 \| \e_{\bX}\{ \chi_1^S(\btheta) \}\|_2 \leq  \delta( \tau, r) r.
\end{align}

Turning to the stochastic component $\bD_0^{-1} \bU(\btheta) = \chi_1^S(\btheta) - \e_{\bX}\{\chi_1^S(\btheta)\}$,
we aim to bound $\max_{\btheta \in \Theta_0(r)} \| \bD_0^{-1} \bU(\btheta) \|_2$, which can be written as
\begin{align}
	\max_{\btheta \in \Theta_0(r), \| \bu \|_2\leq 1} \bu^\T  \bD_0^{-1} \bU(\btheta)   = r^{-1} \max_{\bu, \btheta \in \Theta_0(r)}  \bv^\T \bU(\btheta).   \label{equiv.form}
\end{align}
Note that $\{\bv^\T \bU(\btheta): \bv, \btheta\in \bbr^s\}$ is a bivariate process indexed by $(\bv^\T, \btheta^\T)^\T \in \bbr^{2s}$. Define
\begin{align}
 \bar{\btheta} = (\bv^\T, \btheta^\T)^\T \in \bbr^{2s} , \quad  \bar{\bD}_0 =  \left(
\begin{array}{cc}
\bD_0 &   \mo  \\
 \mo  & \bD_0
\end{array}
\right)   \in \bbr^{(2s)\times (2s)},  \nn  \\
  \bar{\bU}(\bar \btheta) = \bv^\T \bU(\btheta),  \quad  \bar{\Theta}_0(r) = \{ \bar{\btheta} \in \bbr^{2s} :  \| \bar{\bD}_0  \bar{\btheta} \|_2 \leq r \}. \nn
\end{align}
In this notation, from \eqref{equiv.form} and the identity $\bar{\bD}_0 \bar{\btheta} = \bD_0 \bv + \bD_0 \btheta$, it is easy to see that
\begin{align}
	\max_{\btheta \in \Theta_0(r)} \| \bD_0^{-1} \bU(\btheta) \|_2   \leq  r^{-1}  \max_{\bar{\btheta }  \in \bar{\Theta}_0(2r)  }  \bar{\bU}(\bar \btheta).   \label{bilinear.form}
\end{align}

Recall that $\nabla \bzeta^S(\btheta) - \nabla \bzeta^S(\mo)   = -2 \sn \{ I(Y_i\leq  \bX_{iS}^\T \btheta ) -  I(Y_i\leq  0 ) + 1/2 - F_\varepsilon( \bX_{iS}^\T \btheta ) \} \bX_{iS}$, where for $i=1,\ldots , n$, $I(Y_i\leq  \bX_{iS}^\T \btheta ) -  I(Y_i\leq  0 ) + 1/2 - F_\varepsilon(  \bX_{iS}^\T \btheta ) $ is equal to
\begin{align}
 \begin{cases}
 	  I(  0< Y_i\leq  \bX_{iS}^\T \btheta ) - \P_{\bX}( 0< Y_i \leq  \bX_{iS}^\T \btheta )   & \ \ \mbox{ if }   \bX_{iS}^\T \btheta \geq 0, \\
 	- I(  \bX_{iS}^\T \btheta <  Y_i \leq 0) + \P_{\bX}(   \bX_{iS}^\T \btheta < Y_i \leq 0 )   & \ \ \mbox{ if }   \bX_{iS}^\T \btheta < 0.
 \end{cases} \nn
\end{align}
For $\btheta \in \bbr^s$, define random variables $\varepsilon_{i,\btheta}  = I(0<Y_i \leq  \bX_{iS}^\T \btheta )-I(  \bX_{iS}^\T \btheta <Y_i\leq 0)$ satisfying
\begin{itemize}
\item[(i)] conditional on $ \bX_{iS}^\T \btheta \geq 0$, $\varepsilon_{i,\btheta} =1$ with probability $P_{i,\btheta} -1/2$ and $\varepsilon_{i,\btheta}  = 0$ with probability $3/2 - P_{i,\btheta} $;
\item[(ii)] conditional on $ \bX_{iS}^\T \btheta < 0$, $\varepsilon_{i,\btheta} = -1$ with probability $1/2 - P_{i,\btheta} $ and $\varepsilon_{i,\btheta}  = 0$ with probability $1/2 + P_{i,\btheta} $,
\end{itemize}
where $P_{i,\btheta} = F_\varepsilon( \bX_{iS}^\T \btheta )$. In this notation, $\nabla \bzeta^S(\btheta) - \nabla \bzeta^S(\mo) = -2 \sn ( {\rm Id} -\e_{\bX})\varepsilon_{i,\btheta} \bX_{iS}$. For every $\lambda  \in \bbr$ and $\bu \in \bbr^s$, we have
\begin{align}
 & 	\e_{\bX} \exp[ \lambda \bu^\T \{\nabla \bzeta^S(\btheta) - \nabla \bzeta^S(\mo) \} ] \nn \\
 & = \prod_{i=1}^n \Big[ \e_{\bX} \{ e^{  -2 \lambda \bu^\T \bX_{iS}  (I-\e_{\bX})\varepsilon_{i,\btheta} } \} I(   \bX_{iS}^\T \btheta \geq 0) +\e_{\bX} \{ e^{  -2 \lambda \bu^\T \bX_{iS}  (I-\e_{\bX} )\varepsilon_{i,\btheta} } \}  I( \bX_{iS}^\T \btheta < 0)  \Big] \nn \\
 & = \prod_{i=1}^n  \Big[  \big\{ e^{-2\lambda \bu^\T \bX_{iS} (3/2-P_{i,\btheta})} (P_{i,\btheta}-1/2) + e^{2\lambda \bu^\T \bX_{iS}(P_{i,\btheta} -1/2)}(3/2-P_{i,\btheta}) \big\} I( \bX_{iS}^\T \btheta \geq 0)  \nn  \\
 & \qquad \qquad +  \big\{ e^{2\lambda \bu^\T \bX_{iS} (1/2+ P_{i,\btheta})} (1/2 - P_{i,\btheta} ) + e^{2\lambda \bu^\T \bX_{iS}(P_{i,\btheta} -1/2)}(1/2+ P_{i,\btheta})  \big\} I( \bX_{iS}^\T \btheta < 0) \Big]. \nn
\end{align}
Further, using the inequalities $|e^u- 1 - u | \leq \frac{1}{2} u^2  e^{u \vee 0}$ and $1+ u \leq e^u$ which hold for all $u \in \bbr$, the last term above can be bounded by
\begin{align}
 & \prod_{i=1}^n \Big[  \big\{  1 + 2\lambda^2(\bu^\T \bX_{iS})^2 (P_{i,\btheta}-1/2)(3/2-P_{i,\btheta})e^{2\lambda |\bu^\T \bX_{iS}|} \big\}   I( \bX_{iS}^\T \btheta \geq 0)    \nn \\
 & \qquad \qquad +  \big\{  1 + 2\lambda^2(\bu^\T \bX_{iS})^2 (1/2- P_{i,\btheta} )(1/2 + P_{i,\btheta}) e^{2\lambda |\bu^\T \bX_{iS}|} \big\}  I( \bX_{iS}^\T \btheta < 0)  \Big]  \nn \\
 & \leq \prod_{i=1}^n \big\{  1+ 2 \lambda^2  (\bu^\T \bX_{iS})^2 |P_{i,\btheta} -1/2| e^{2\lambda |\bu^\T \bX_{iS}|} \big\}  \nn \\
 & \leq  \prod_{i=1}^n \exp\big\{ 2 \lambda^2  (\bu^\T \bX_{iS})^2 |P_{i,\btheta} -1/2| e^{2\lambda |\bu^\T \bX_{iS}|} \big\} .  \nn
\end{align}
Consequently, for every $\bar{\btheta}  = (\bv^\T, \btheta^\T)^\T \in \bar{\Theta}_0(2r)$,
\begin{align}
	& \log \e_{\bX} \exp\bigg\{  \lambda \frac{\bar{\bU}(\bar{\btheta}) - \bar{\bU}(\mo) }{ \| \bar{\bD}_0 \bar{\btheta} \|_2 } \bigg\} = \log \e_{\bX} \exp\bigg\{  \lambda \frac{ \bv^\T \{ \bzeta^S(\btheta) - \bzeta^S(\mo)\} }{  \| \bar{\bD}_0 \bar{\btheta} \|_2  }  \bigg\}  \nn \\
	& \leq  \frac{2\lambda^2}{\| \bD_0 \bv \|_2^2 + \| \bD_0 \btheta\|_2^2 }\sn (\bv^\T \bX_{iS})^2 | P_{i,\btheta}-1/2| \exp\bigg(  \frac{2\lambda |\bv^\T \bX_{iS}|}{\| \bar{\bD}_0 \bar{\btheta} \|_2} \bigg) .  \label{mgf.bound.1}
\end{align}
On the event $\mathcal{E}_0(\tau )$ for some $\tau >0$, we have $|P_{i,\btheta} - 1/2 | \leq 2A_2  \{2n f_\varepsilon(0)\}^{-1/2}  \tau^{1/2} r $ and $ |\bv^\T \bX_{iS}|  \leq  \| \bD_0 \bv \|_2 \| \bD_0^{-1} \bX_{iS} \|_2 \leq \| \bD_0 \bv \|_2 \{ 2n f_\varepsilon(0) \}^{-1/2} \tau^{1/2}$. Together with \eqref{mgf.bound.1}, this yields that for all $|\lambda | \leq \{ 2n f_\varepsilon(0) / \tau \}^{1/2}$,
\begin{align}
   \log \e_{\bX} \exp\bigg\{  \lambda \frac{\bar{\bU}(\bar{\btheta}) - \bar{\bU}(\mo) }{ \| \bar{\bD}_0 \bar{\btheta} \|_2 }  \bigg\}  \leq  \frac{\lambda^2}{2}   \frac{4 e^2 A_2\, r}{f_\varepsilon(0)}  \sqrt{  \frac{\tau}{2n f_\varepsilon(0)}  }.  \label{mgf.bound.2}
\end{align}
In view of \eqref{mgf.bound.2}, define
\begin{align}
	w_0(\tau) =  2 e  \sqrt{  \frac{A_2  \, r_0 }{f_\varepsilon(0)} }   \bigg\{ \frac{\tau }{2n f_\varepsilon(0)} \bigg\}^{1/4}   \label{def.w}
\end{align}
for some $r_0>0$ to be specified (see \eqref{def.r0} below), such that for any $\bar{\btheta}  = (\bv^\T, \btheta^\T)^\T \in \bar{\Theta}_0(2r)$ with $0\leq r\leq r_0$,
\begin{align}
   \e_{\bX} \exp\bigg\{ \frac{\lambda}{w_0(\tau)} \frac{\bar{\bU}(\bar{\btheta}) - \bar{\bU}(\mo) }{    \| \bar{\bD}_0 \bar{\btheta} \|_2 }  \bigg\} \leq  \exp(  \lambda^2 /2 )   \label{mgf.bound.3}
\end{align}
holds almost surely on $\mathcal{E}_0(\tau )$ for all
\begin{align}
	|\lambda| \leq    2 e   \sqrt{  \frac{A_2  \, r_0 }{f_\varepsilon(0)} }  \bigg\{ \frac{2n f_\varepsilon(0)}{\tau } \bigg\}^{1/4} := g_0(\tau)  . \label{def.g}
\end{align}
By \eqref{mgf.bound.3}, it follows from Corollary~2.2 in the supplement of \cite{S2012} and \eqref{bilinear.form} that, for any $\tau >0$, $0\leq r\leq r_0$ and $0<t \leq  \frac{1}{2} g_0^2(\tau)   - 2s$,
\begin{align}  \label{stochastic.tail.bound}
	\P_{\bX}  \bigg\{  \max_{ \btheta  \in  {\Theta}_0( r ) } \| \bD_0^{-1} \bU(\btheta)  \|_2    \geq 6 w_0(\tau)   ( 2t + 4s )^{1/2}  \bigg\} \leq e^{-t}
\end{align}
holds almost surely on $\mathcal{E}_0(\tau )$, where $g_0$ is given at \eqref{def.g}.

Combining \eqref{mean.chi-theta} and \eqref{stochastic.tail.bound} we obtain that for any $\tau >0$, $0\leq r\leq r_0 \leq \{ 2n f_\varepsilon(0)  /\tau \}^{1/2}$ and $0< t \leq  \frac{1}{2} g_0^2(\tau)  - 2s $,
\begin{align}
 \P_{\bX} \bigg\{  \max_{\btheta \in \Theta_0(r) } \| \chi_1^S(\btheta)  \|_2  \geq \delta(\tau, r ) r + 6 w_0(\tau) (2t + 4s )^{1/2} \bigg\} \leq e^{-t}  \label{chi1.tail.bound}
\end{align}
almost surely on $\mathcal{E}_0(\tau )$. For a given triplet $(\tau , r, t)$, define the event
\begin{align}
	\Omega^S_0(\tau , r, t) = \bigg\{  \max_{\btheta \in \Theta_0(r) } \| \chi^S_1(\btheta)  \|_2  \leq \delta( \tau, r ) r + 6 w_0(\tau) (2t + 4s )^{1/2} \bigg\}. \label{event.Omega}
\end{align}

\medskip
\noindent
{\it \textbf{Step~2: Fisher approximation}}. By Lemma~\ref{LAD.concentration},
\begin{align}
	&  \max_{S\subseteq [p]: |S|=s} \| \bD_0 \hat{\btheta}_S \|_2 \nn \\
	& = \{ 2 n f_\varepsilon(0) \}^{1/2} \max_{S\subseteq [p]: |S|=s} \| \hat \bSigma_{SS}^{1/2}  \hat{\btheta}_S \|_2 \leq   C_1 \, a_2^{-1} \{ 2 f_\varepsilon(0) s\log(pn) \}^{1/2} := r_0 \label{def.r0}
\end{align}
holds with probability at least $1-c_1 n^{-1} - c_2 n^{1-\kappa}$. Moreover, since $\hat{\btheta}_S$ maximizes ${\cal L}_n^S(\btheta)$ over $\btheta\in \bbr^s$ for each $s$-subset $S\subseteq [p]$, we have $\nabla {\cal L}_n^S(\hat{\btheta}_S) = \mo$ and $\chi_1^S(\hat{\btheta}) =  \bD_0 \hat{\btheta}_S - \hat{\bxi}_S$. This, together with \eqref{event.Omega} implies that on the event $\{ \hat{\btheta}_S \in \Theta_0(r_0) \} \cap \Omega_0^S(\tau , r_0,t)$,
\begin{align}
	   \| \bD_0 \hat{\btheta}_S - \hat{\bxi}_S \|_2 \leq \delta(\tau, r_0) r_0 + 6 w_0(\tau) ( 2t + 4s )^{1/2} \label{Fisher.approxi}
\end{align}
whenever $n \geq \{2f_\varepsilon(0)\}^{-1}   \tau r_0^2 $.

\medskip
\noindent
{\it \textbf{Step~3: Wilks approximation}}. For $ \btheta_1, \btheta_2 \in \Theta_0(r)$, define
\begin{align}
\chi_2^S(\btheta_1, \btheta_2) & = {\cal L}_n^S(\btheta ) - {\cal L}_n^S( {\btheta}_2) - ( \btheta_1 -  \btheta_2 )^\T \nabla {\cal L}_n^S( {\btheta}_2) + \frac{1}{2} \| \bD_0 ( \btheta_1 -  {\btheta}_2 ) \|_2^2 . \label{def.chi2-2}
\end{align}
Noting that $\nabla_{\btheta_1} \chi_2^S(\btheta_1, \btheta_2)  = \nabla {\cal L}_n^S(\btheta_1)- \nabla {\cal L}_n^S( \btheta_2 ) + \bD_0^2 (\btheta_1 - \btheta_2 ) = \bD_0 \{  \chi_1^S( \btheta_1 )  - \chi_1^S( \btheta_2 ) \}$, we have
\begin{align}
	  |\chi_2^S(\btheta_1, \btheta_2) |  = |\chi_2^S(\btheta_1 , \btheta_2 ) - \chi_2^S(\btheta_2, \btheta_2 )|  \leq  2\| \bD_0 ( \btheta_1 - \btheta_2 )\|_2  \max_{\bu \in \Theta_0(r) } \| \chi_1^S(\bu) \|_2, \label{chi2.ubd}
\end{align}
where $\wt \btheta = \lambda \btheta $ for some $0\leq \lambda \leq 1$. Let $r_0 >0$ be as in \eqref{def.r0}. Then, it follows from \eqref{chi2.ubd} that on $\Omega^S_0(\tau , r_0, t)$ with $n \geq \{2f_\varepsilon(0)\}^{-1}   \tau r_0^2$,
\begin{align}
	\max_{\btheta_1, \btheta_2 \in \Theta_0(r_0)} \frac{|\chi_2^S(\btheta_1, \btheta_2)|}{\| \bD_0( \btheta_1 - \btheta_2 ) \|_2 } \leq 2\delta(\tau, r_0 ) r_0 + 12 w_0(\tau) ( 2t + 4s )^{1/2} . \nn
\end{align}

In view of \eqref{def.chi2-2}, ${\cal L}_n^S(\hat{\btheta}_S) - {\cal L}_n^S(\mo) - \frac{1}{2} \| \bD_0 \hat{\btheta}_S \|_2^2  = - \chi_2^S(\mo, \hat{\btheta}_S)$. Therefore, on the event $\{ \hat{\btheta}_S \in \Theta_0(r_0) \} \cap \Omega^S_0(\tau , r_0,t)$ we have
\begin{align}
	& \Big|   [ 2\{ {\cal L}_n^S(\hat{\btheta}_S) - {\cal L}_n^S(\mo) \} ]^{1/2}  - \| \bD_0 \hat{\btheta}_S \|_2  \Big| \nn \\
	& \leq \frac{|2\{ {\cal L}_n^S(\hat{\btheta}_S) - {\cal L}_n^S(\mo) \} -  \|  \bD_0 \hat{\btheta}_S\|_2^2  |}{\| \bD_0 \hat{\btheta}_S \|_2} \leq \frac{2|\chi_2^S(\mo, \hat{\btheta}_S)|}{ \| \bD_0 \hat{\btheta}_S \|_2 } \leq 4  \big\{ \delta(\tau, r_0 ) r_0 + 6 w_0(\tau)  ( 2t + 4s )^{1/2}  \big\} ,  \nn
\end{align}
provided that $n \geq \{2f_\varepsilon(0)\}^{-1}  \tau r_0^2 $. Together with \eqref{Fisher.approxi}, this implies that conditional on the event $\cap_{S\subseteq [p]: |S| =s } \{ \hat{\btheta}_S \in \Theta_0(r_0) \} \cap \Omega^S_0(\tau , r_0,t)$,
\begin{align}
 \max_{S\subseteq [p] : |S| = s }  \Big|  [ 2\{ {\cal L}_n^S(\hat{\btheta}_S) - {\cal L}_n^S(\mo) \}   ]^{1/2}  - \| \hat{\bxi}_S \|_2  \Big| \leq  5 \big\{ \delta(\tau, r_0 ) r_0 + 6 w_0(\tau) ( 2t + 4s )^{1/2} \big\}   \label{sqrt.approxi}
\end{align}
whenever $n \geq \{2f_\varepsilon(0)\}^{-1}  r_0^2 \tau $, where $\delta(\tau, r)$, $r_0$ and $w_0(\tau)$ are as in \eqref{def.delta-r}, \eqref{def.r0} and \eqref{def.w}.

Finally,  taking $\tau =\tau_0  \asymp  \lambda^{-1}_{\min}(s)  s\log(pn)$ as in \eqref{uniform.Wilks.approxi} and setting $t=s\log \frac{ep}{s} + \log n$ in the concentration bound \eqref{chi1.tail.bound} prove \eqref{wilks.approx.2} using Boole's inequality. \qed

\end{document}